\documentclass{scrartcl}
\usepackage{amsmath}
\usepackage{amsfonts}
\usepackage{amssymb}
\usepackage{dsfont}
\usepackage[T1]{fontenc}
\usepackage[latin1]{inputenc}
\usepackage{epsfig}
\usepackage{graphicx}

\usepackage[boxed, lined]{algorithm2e}
\usepackage{float}
\usepackage{comment}
\usepackage{mathtools}

\usepackage{bm}

\newcommand{\subalign}[1]{%
   \vcenter{%
     \Let@ \restore@math@cr \default@tag
     \baselineskip\fontdimen10 \scriptfont\tw@
     \advance\baselineskip\fontdimen12 \scriptfont\tw@
     \lineskip\thr@@\fontdimen8 \scriptfont\thr@@
     \lineskiplimit\lineskip
     \ialign{\hfil$\m@th\scriptstyle##$&$\m@th\scriptstyle{}##$\hfil\crcr
       #1\crcr
     }%
   }%
}

\newcommand{\bw}{\mathbf{w}}

\newcommand{\D}{{\mathcal{D}}}

\newcommand{\W}{{\mathcal{W}}}

\newcommand{\Nu}{{\mathcal{N}}}

\newcommand{\N}{\mathbb{N}}

\newcommand{\R}{\mathbb{R}}

\newcommand{\Rd}{\mathbb{R}^d}

\newcommand{\beq}{\begin{eqnarray*}}

\newcommand{\eeq}{\end{eqnarray*}}

\newcommand{\beqm}{\begin{eqnarray}}

\newcommand{\eeqm}{\end{eqnarray}}

\newtheorem{theorem}{Theorem}
\newtheorem{corollary}{Corollary}
\newtheorem{lemma}{Lemma}
\newtheorem{definition}{Definition}

\DeclareOldFontCommand{\bf}{\normalfont\bfseries}{\mathbf}
\DeclareOldFontCommand{\it}{\normalfont\itshape}{\mathit}

\newcommand{\EXP}{{\mathbf E}}
\newcommand{\PROB}{{\mathbf P}}

\renewcommand{\P}{{\cal P}}
\renewcommand{\bf}{\normalfont \bfseries}
\renewcommand{\it}{\normalfont \itshape}

\allowdisplaybreaks
\begin{document}
\renewcommand{\thefootnote}{\fnsymbol{footnote}}
\newcommand{\F}{{\cal F}}
\newcommand{\Sp}{{\cal S}}
\newcommand{\G}{{\cal G}}
\newcommand{\HH}{{\cal H}}

\begin{center}

  {\LARGE \bf
    Statistically guided deep learning
  }
\footnote{
Running title: {\it Statistically deep learning}}
\vspace{0.5cm}

Michael Kohler $^{1}$ and
Adam Krzy\.zak$^{2,}$\footnote{Corresponding author. Tel:
  +1-514-848-2424 ext. 3007, Fax:+1-514-848-2830}

{\it $^1$
Fachbereich Mathematik, Technische Universit\"at Darmstadt,
Schlossgartenstr. 7, 64289 Darmstadt, Germany,
email: kohler@mathematik.tu-darmstadt.de
}

{\it $^2$ Department of Computer Science and Software Engineering, Concordia University, 1455 De Maisonneuve Blvd. West, Montreal, Quebec, Canada H3G 1M8, email: krzyzak@cs.concordia.ca}

\end{center}
\vspace{0.5cm}

\begin{center}
April 11, 2025
\end{center}
\vspace{0.5cm}

\noindent
    {\bf Abstract}\\
    We present a theoretically well-founded deep learning algorithm
    for nonparametric regression. It uses over-parametrized
    deep neural networks with logistic activation function, which
    are fitted to the given data via gradient descent.
    We propose a special topology of these networks, a special
    random initialization of the weights, and a data-dependent choice
    of the learning rate and the number of gradient descent steps.
    We prove a theoretical bound on the expected
    $L_2$ error of this estimate, and illustrate its finite sample size
    performance by applying it to simulated data.

    Our results show that a theoretical analysis of deep learning
    which takes into account simultaneously optimization, generalization
    and approximation can result in a new deep
    learning estimate which has an improved finite
    sample performance.
    
    \vspace*{0.2cm}


\vspace*{0.2cm}


\section{Introduction}
\label{se1}

\subsection{Scope of this paper}
\label{se1sub1}
Due to its tremendous success in applications, e.g., in
image classification
(cf., e.g., Krizhevsky, Sutskever and Hinton  (2012)),
in language recognition (cf., e.g.,  Kim (2014))
in machine translation (cf., e.g., Wu et al. (2016))
or in mastering of games (cf., e.g., Silver et al. (2017)),
deep learning is currently changing the world. This big success
of deep learning in the past relies on two things: the massive
increase of computing power and availability of the huge data
sets. However, it seems that both cannot be much more increased:
Firstly, there is already a shortage of computer chips for deep learning,
and also the increasing electricity demand of the computers used
for computing the deep learning estimates seems problematic.
And secondly, e.g. for large language models, all available
text data has been already used for the training, so it is not clear
how the size of the used data sets can be further increased. 

But there remains one different
approach to improve the deep learning estimates: one can try to improve
the used estimation methods. In the past new methods have been mainly
constructed by trial and error, and not based on a rigorous
theoretical analysis.
In this paper we investigate whether a theoretical approach 
succeeds in improving deep learning estimates.

\subsection{Nonparametric regression}
\label{se1sub2}
We study deep learning estimates in the context of nonparametric
regression. Here
$(X,Y)$, $(X_1,Y_1)$, \dots, $(X_n,Y_n)$
are independent and identically $\R^d \times
\R$--valued
random vectors with $\EXP Y^2 < \infty$, and given the data set
\begin{equation}
\label{inteq1}
\D_n=\{(X_1,Y_1), \dots, (X_n,Y_n)\}
\end{equation}
the task is to estimate the so--called regression function
\[
m:\R^d \rightarrow \R, \quad m(x)=\EXP\{Y|X=x\}.
\]
More precisely, the goal is to
construct an estimate
\[
m_n(\cdot)=m_n(\cdot, \D_n): \R^d \rightarrow \R
\]
such that the so-called $L_2$ error 
\[
\int | m_n(x)-m(x)|^2 \PROB_X(dx)
\]
is close to zero.

A detailed introduction to nonparametric regression, its estimates
and known theoretical results can be found, e.g., in  Gy\"orfi et al. (2002).

\subsection{Least squares estimates estimates}
\label{se1sub3}
Since
\[
\EXP\{ |m_n(X)-Y|^2 | \D_n\}
=
\EXP\{ |m(X)-Y|^2 \}
+
\int | m_n(x)-m(x)|^2 \PROB_X(dx)
\]
(cf., e.g., Chapter 1 in Gy\"orfi et al. (2002)),
the aim of minimizing the $L_2$ error
means that one
wants
to find an estimate such that its so--called $L_2$ risk (or mean
squared
prediction error)
\begin{equation}
\label{inteq2}
\EXP\{ |m_n(X)-Y|^2 | \D_n\}
\end{equation}
is close to the optimal value $\EXP\{ |m(X)-Y|^2 \}$.

This way of considering the estimation task immediately suggest a way
of solving it: One can try to use the given data (\ref{inteq1})
to estimate the $L_2$ risk (\ref{inteq2}) by the so--called empirical
$L_2$ risk
\begin{equation}
\label{inteq3}
\frac{1}{n} \sum_{i=1}^n | m_n(X_i) - Y_i|^2
\end{equation}
and can try  to minimize (\ref{inteq3}) over some space of functions.
This leads to so--called least squares estimates
\begin{equation}
  \label{inteq4}
  m_n(\cdot) = \arg \min_{f \in \F_n}
\frac{1}{n} \sum_{i=1}^n | f(X_i) - Y_i|^2
\end{equation}
which depend on spaces $\F_n$ of functions $f:\R^d \rightarrow \R$.
Here the right choice of the function space is crucial, since it must
be on the one hand so rich that functions in it are able to approximate
the (unknown) regression function well, and on the other hand
it should be such  that the empirical $L_2$ risk
of  the function which
minimizes the empirical $L_2$ risk
is close
to its expectation. Usually the latter is shown by
showing that the maximal deviation between the $L_2$ risk and
the empirical $L_2$ risk
on the function space
is small, which holds if the function
space is not too complex.

\subsection{Neural networks}
\label{se1sub4}
For neural network estimates one considers in this context spaces
of neural networks.
In their simplest form of fully connected feedfoward neural networks they
are defined as follows: One chooses
an activation function
$\sigma: \R \rightarrow \R$, e.g.,
\begin{equation}
\label{inteq5}
\sigma(x)=\max \{x,0\}
\end{equation}
(so-called ReLU-activation function) or
\begin{equation}
\label{inteq6}
\sigma(x)= \frac{1}{1+e^{-x}}
\end{equation}
(so-called logistic squasher), and selects the number $L
\in \N$ of hidden layers of the network and the numbers $k_s \in \N$
of neurons in the $s$-th hidden layer $(s \in \{1, \dots, L\})$.
Then the feedforward neural network $f_\bw$ with $L$ hidden layers,
$k_s$ neurons in layer $s \in \{1, \dots, L\}$ and with weight vector
$\bw=(w_{i,j}^{(l)})_{l,i,j}$ is the function $f_\bw:\R^d \rightarrow
\R$
defined by
\begin{equation}
  \label{inteq7}
  f_\bw(x)=
\sum_{j\in \{1, \dots,k_L\}}
w_{1,j}^{(L)} \cdot f_j^{(L)}(x),
\end{equation}
where
\begin{equation}
\label{inteq8}
f_i^{(s)}(x)
=
\sigma \left(
\sum_{
j \in \{1, \dots,k_{s-1}\}}
w_{i,j}^{(s-1)} \cdot f_j^{(s-1)}(x)
+ w_{i,0}^{(s-1)}
\right)
\quad \mbox{for }
s \in \{2, \dots, L\} \mbox{ and } i>0
\end{equation}
and
\begin{equation}
\label{inteq9}
f_i^{(1)}(x)
=
\sigma \left(
\sum_{
j \in \{1, \dots, d\}
}
w_{i,j}^{(0)} \cdot x^{(j)} + w_{i,0}^{(0)} 
\right)
\quad \mbox{for }
i>0.
\end{equation}
Here $w_{i,j}^{(s-1)}$ is the weight between neuron $j$ in layer
$s-1$ and neuron $i$ in layer $s$. And  $w_{i,0}^{(s-1)}$
is the bias in the computation of the output of neuron $i$ in layer
$s$.

The idea is then to fix the activation function, the numer of layers $L \in \N$,
the number $k_s \in \N$ of neurons in layer $s \in \{1, \dots, L\}$, and
to choose the weight vector $\bw$ by minimizing
the empirical $L_2$ risk
\begin{equation}
\label{inteq10}
F_n(\bw) = \frac{1}{n} \sum_{i=1}^n | f_\bw(X_i)-Y_i|^2
\end{equation}
of $f_\bw$
with respect to $\bw$. 

Usually, the activation function $\sigma$ is highly nonlinear, and
therefore
$f_{\bw}(X_i)$ and also $F_n(\bw)$ depend nonlinearly on $\bw$. Due
to this it is not clear how one can minimize (\ref{inteq10}) with
respect
to the weight vector $\bw$.

\subsection{Computation of neural network regression estimates}
\label{se1sub5}
Minimizing of the empirical $L_2$ risk with respect to a class of neural
networks is done in practice by using gradient descent
(or one of its variants like stochastic gradient descent):
One chooses a random
starting vector $\bw^{(0)}$ for the weights and  computes
$t_n \in \N$ gradient descent steps
\begin{equation}
\label{inteq11}
\bw^{(t)}
=
\bw^{(t-1)}
-
\lambda_n \cdot \nabla_{\bw} F_n(\bw^{(t-1)})
\quad
(t=1, \dots, t_n)
\end{equation}
with stepsize $\lambda_n>0$. The estimate is then defined by
\[
m_n(x)=f_{\bw^{(t_n)}}(x).
\]

\subsection{Difficulty in the application of deep neural networks}
\label{se1sub6}

The above definition of the neural network regression estimates requires
decisions about the class of neurals networks which we fit to the data, the
choice of the starting vector, the choice of the number of gradient descent
steps and the choice of the stepsize.

If we consider for simplicity just fully connected neural networks with
$L$ layers and $r$ neurons per layer
(i.e., we set $k_s=r$ for $s=1, \dots, L$), the logistic activation function,
and the famous ADAM rule for the choice of the stepsize, then the
remaining question is  how to choose the starting vector, and how to choose
the number of gradient descent steps. For the choice of the starting vector
popular algorithms in the literature are the \linebreak GlorotNormal-,
the GlorotUniform-, the
HeNormal- or the HeUniform-rule (cf., e.g.,
Chapter 8 in Goodfellow, Bengio and Courville (2016)), where the initial weights
are chosen independently from the normal distributions or uniform distributions.

        	\begin{figure}[h!]
	  \begin{tabular}{cc}
	    \includegraphics[width=7cm]{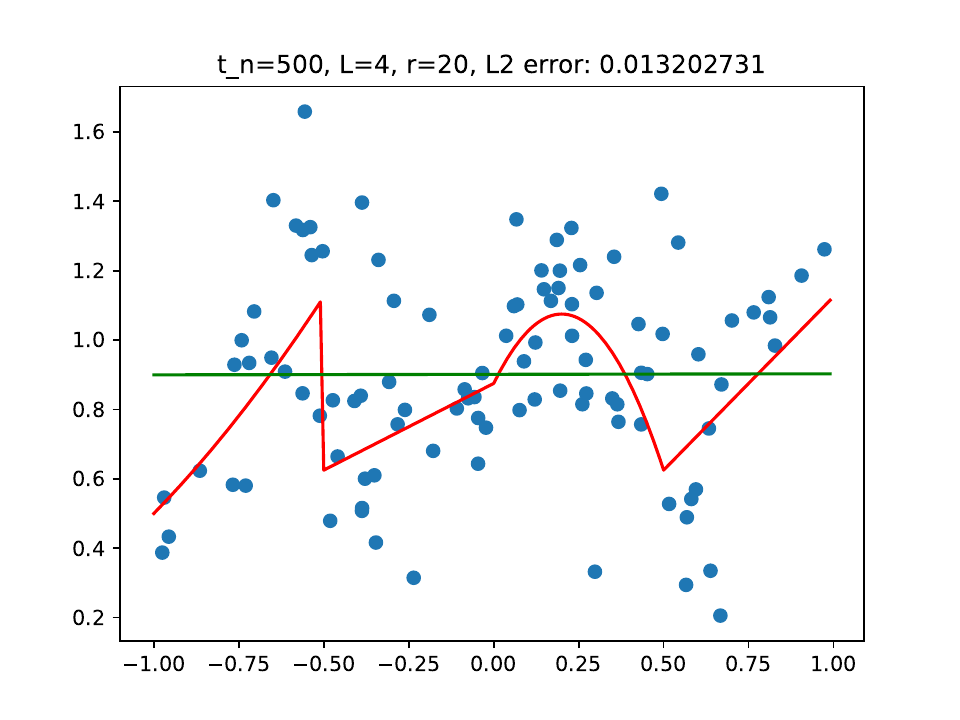}
            &
            		\includegraphics[width=7cm]{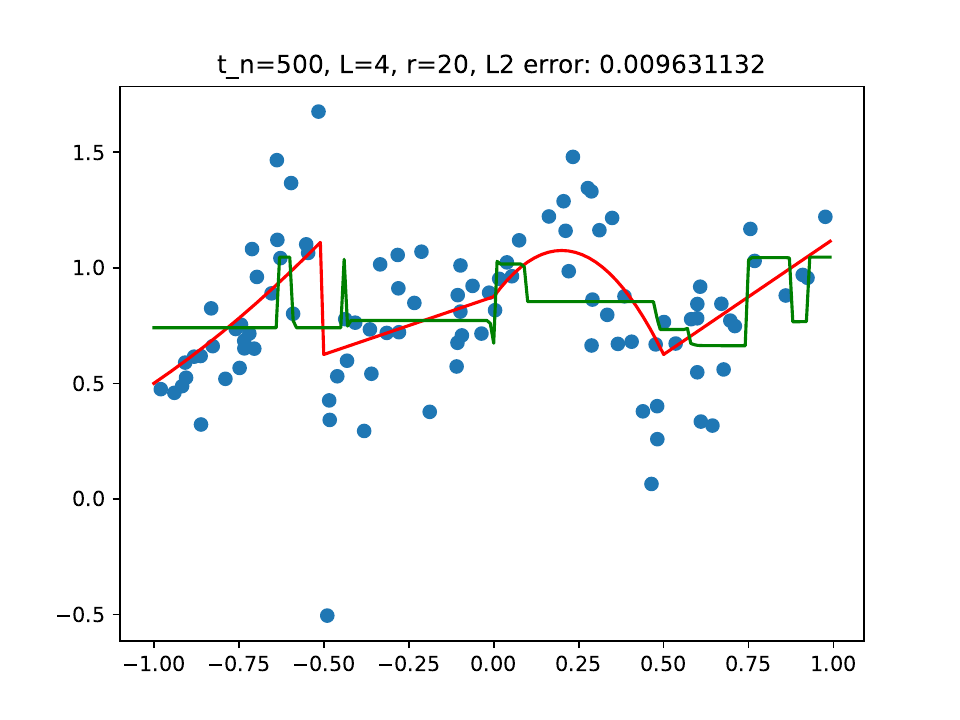}
	                \\
\includegraphics[width=7cm]{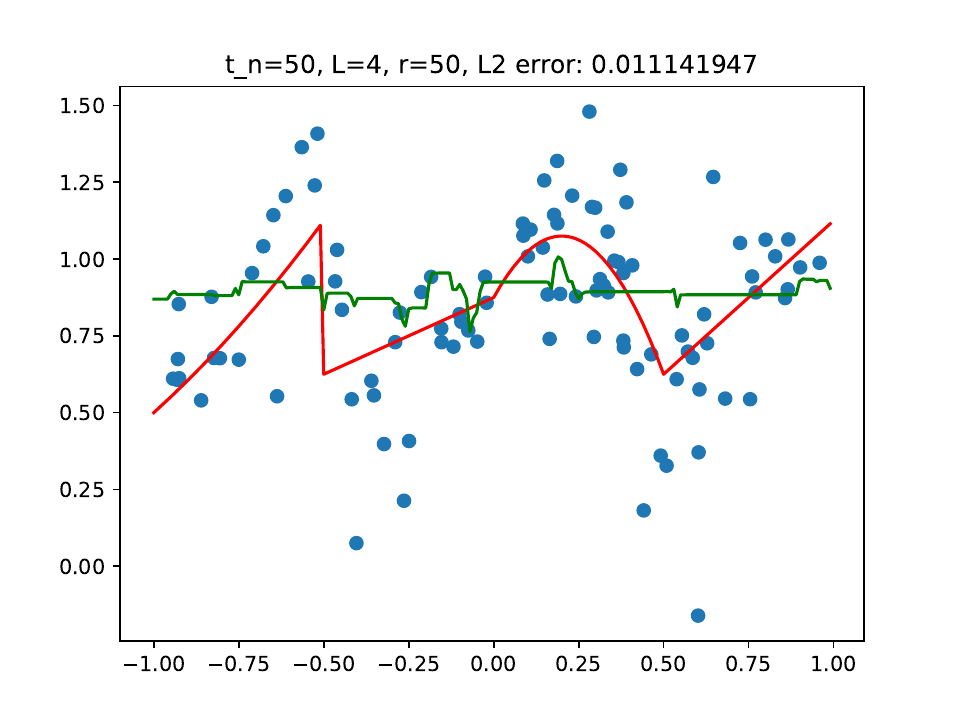}
&
\includegraphics[width=7cm]{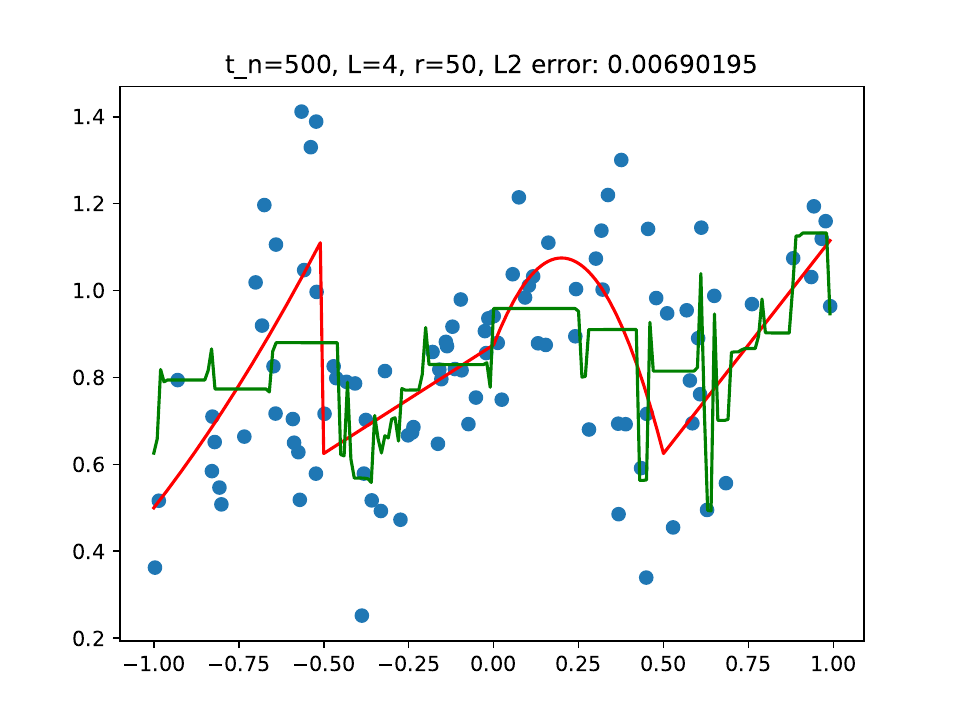}
        \end{tabular}
	  \caption{Neural network estimate with various initialization schemes,
            various topologies and various choices of the stepsize applied to
            the univariate regression problem with sample size $n=100$. \label{fig_int_1}}
	\end{figure}

In the upper right panel in Figure \ref{fig_int_1} we apply this for a neural network with $L=4$, $r=20$ and the GlorotNormal-rule for initialization
to an univariate regression problem (which is described in detail in Section \ref{se4}), which leads to a constant estimate (in green) which does not approximate the regression function (in red) well. The same effect happens with
the GlorotUniform-, the
HeNormal- or the HeUniform-rule.
The picture drastically changes if we use the uniform distribution on an interval
$[-A,A]$ for the weights on the input level, the uniform distribution
on an interval $[-B,B]$ for all inner weights, set all weights on the output level initially to zero, and choose  a large value for $A$ and  a
moderate value for $B$. For $A=1000$, $B=20$ and three different values for $(L,r,t_n)$
the estimates are then shown in the upper right, the lower left, and the lower
right panel in Figure \ref{fig_int_1}, respectively.

This shows that the performance of the neural network estimate crucially
depends on the chosen parameters.
As mentioned on page 293 in Goodfellow, Bengio and Courville (2016)
``designing improved initialization strategies is a difficult
task because neural network optimization is not yet well understood.''
Furthermore, it is mentioned there that ``our understanding of how the
initial point affects generalization is especially primitive, offering
little or no guidance for how to select the initial point''.
It should be added that  the same problem also occurs in
connection with the number
of gradient descent steps for the ADAM rule, or more generally
with the stepsize choices during gradient descent
and the number of gradient descent steps.

In this article we consider simultaneously optimization, generalization
and approximation of neural networks and use this to
propose a theoretically motivated way of choosing the parameters
of the neural network estimates.

\subsection{A theoretical approach to deep learning}
\label{se1sub7}

In practice usually over-parametrized deep neural networks
are used, where the number of weights is much larger than the
sample size $n$, so one fits a function to the data which has
much more free parameters (i.e., weights) than there are data
points.

There are three main theoretical questions in this context: Firstly, why does the
resulting estimate optimize well, i.e., why
is gradient descent able to achieve small values of the empirical
$L_2$ risk? Secondly, why does the resulting estimate generalize well,
i.e., why is its squared error on new independent data (not contained
in the training data) small? And why does it approximate well,
i.e., why does the sequence of weight vectors considered during gradient
descent contains a weight vector for which the corresponding neural
network approximates the regression function well? Of course, if
we are able to answer these questions then it should also be possible
to say which activation function, which topology  (i.e., number of
layers and number of neurons per layer), which initialization of the weights,
which stepsize and which number of gradient descent steps lead to estimates with 
a small $L_2$ error. So a theoretical understanding of the
above three questions might be used to construct hints for the
choice of the parameters of the estimate in applications.

Kohler (2024) has developed a theory answering these questions,
which applies to over-parametrized deep neural networks with
smooth activation function. It uses the observation that for a
proper choice of $\lambda_n$ and $t_n$ 
the weights computed during gradient descent stay 
in a local neighborhood of the starting value. More precisely,
if $\lambda_n= \frac{1}{L_n}$ and the gradient of $F_n(\bw)$ is 
Lipschitz continuous with Lipschitz constant $L_n$ around the
starting weight vector, i.e., if
\[
\| \nabla_\bw F_n(\bw_1)- \nabla_\bw (F_n(\bw_2) \| \leq L_n \cdot \|
\bw_1 - \bw_2 \|
\]
holds for $\bw_1$ and $\bw_2$ ''close'' to the starting weight vector
$\bw^{(0)}$, and additionally the gradient is suitably bounded in this neighborhood,
then during gradient descent
\[
\|\bw^{(t)}-\bw^{(0)}\|\leq \sqrt{ c_1 \cdot \lambda_n \cdot t}
\]
holds for all $t \in \{1, \dots, t_n\}$ (cf., Lemma A.1 in Braun et
al. (2023)). Since
\[
\|\bw^{(t)}-\bw^{(0)}\|_\infty \leq \|\bw^{(t)}-\bw^{(0)}\|
\]
this implies that if we choose $\lambda_n$ and $t_n$ such that
$\lambda_n \cdot t_n$ is bounded by some constant, then any bounds
which we impose on the absolute value of the weights during the random initialization
enable us to derive bounds on the absolute value of the weights during
gradient descent.

Kohler (2024) uses such bounds to ensure the 
estimates generalize well. This is possible, since the smoothness of the activation
function together with the bounds on the weights enables one
to derive bounds on the derivative of the networks. And using these
bounds one can approximate the 
corresponding set of deep networks by piecewise
polynomials and bound the complexity of the set of deep networks
by a suitable covering number of the set of piecewise polynomials.
In this context Kohler (2024) uses a special topology of the network,
where a huge linear combination of many small networks of fixed
depth $L$ and width $r$ are computed. It turns out that neither the
number $K_n$ of these small networks nor the bounds on the absolute value of
the coefficients in the linear combination have a crucial influence
on the covering number above as long as they grow not faster
than some polynomial in the sample size. In this way it is possible
to define over-parametrized deep neural networks which generalize
well (since during gradient descent they are always contained
in some function space with a finite complexity). Furthermore,
Kohler (2024) uses different bounds $A_n$ and $B_n$ for the
absolute value of the weights in the input layer and the absolute
value of the weights between the hidden layers.
Here $B_n$ is chosen as a large constant, and then $A_n$ is the
main parameter controlling the complexity of the over-parametrized
deep networks chosen by $A_n= c_2 \cdot (\log n) \cdot n^\tau$ for some $\tau \in (0,1)$.

In order to analyze the approximation error, Kohler (2024)
derives an approximation result for the approximation
of a smooth function by networks where the weights are
bounded as above. Here the number of networks and the
size of $A_n$ are related and they control the approximation
error of the deep network.

Furthermore, Kohler (2024) uses a relation between the 
gradient descent applied to the empirical $L_2$ risk
of the deep network and the gradient descent applied
to the empirical $L_2$ risk of the linear Taylor
approximation of the deep network in order to
analyze the gradient descent. This makes it possible
to use techniques which have been developed for analysis of gradient
descent applied to smooth convex functions. 

In Kohler (2024) the
regression function is assumed to be $(p,C)$--smooth in
the following sense.

\begin{definition}
\label{se2de1} 
  Let $p=q+s$ for some $q \in \N_0$ and $0< s \leq 1$.
A function $m:\R^d \rightarrow \R$ is called
$(p,C)$-smooth, if for every $\bm{\alpha}=(\alpha_1, \dots, \alpha_d) \in
\N_0^d$
with $\sum_{j=1}^d \alpha_j = q$ the partial derivative
$\partial^q m/(\partial x_1^{\alpha_1}
\dots
\partial x_d^{\alpha_d}
)$
exists and satisfies
\[
\left|
\frac{
\partial^q m
}{
\partial x_1^{\alpha_1}
\dots
\partial x_d^{\alpha_d}
}
(x)
-
\frac{
\partial^q m
}{
\partial x_1^{\alpha_1}
\dots
\partial x_d^{\alpha_d}
}
(z)
\right|
\leq
C
\cdot
\|x-z\|^s
\]
for all $\bold{x},\bold{z} \in \R^d$, where $\Vert\cdot\Vert$ denotes the Euclidean norm.
\end{definition}

Kohler (2024) considered a neural network topology consisting of
$K_n \in \N$ in parallel computed neural network with
logistic squasher activation function, and with
depth $L$ and width $r$, where
  \[
  \frac{K_n}{n^\kappa} \rightarrow 0 \quad (n \rightarrow \infty)
  \quad \mbox{and} \quad
  \frac{K_n}{n^{4 \cdot r \cdot (r+1) \cdot (L-1) + r \cdot (4d+6) + 6}}
  \rightarrow \infty \quad (n \rightarrow \infty)
  \]
  for some $\kappa>0$ and
\[
L=\lceil \log_2(q+d) \rceil+1 \quad \mbox{and} \quad r=2 \cdot
\lceil (2p+d)^2 \rceil.
\]
The weights are initialized such that all weights of the input
level are uniformly distributed on $[-c_2 \cdot (\log n) \cdot n^{1/(2p+d)},
  c_2 \cdot (\log n) \cdot n^{1/(2p+d)}]$, all inner weights are uniformly distributed
on $[-c_3, c_3]$ and all the output weights are set to zero.
Then
\[
t_n=\left\lceil
c_4 \cdot \frac{K_n^3}{\beta_n}
\right\rceil
\]
gradient descent steps with stepsize
\[
\lambda_n=\frac{c_5}{ n \cdot K_n^3}
\]
are performed. It is shown in Theorem 1 in Kohler (2024)
that a truncated version
of this estimate satisfies
for any $\epsilon>0$
\[
\EXP \int | m_n(x)-m(x)|^2 \PROB_X (dx)
\leq c_6 \cdot n^{- \frac{2p}{2p+d} + \epsilon},
\]
provided $supp(X)$ is compact, $\EXP\left\{
e^{c_7 \cdot Y^2}
\right\}
< \infty$ and the regression functions is $(p,C)$--smooth.
  
The main problem in using this result for defining a neural network
estimate applied to data is that the parameters $K_n$ and $t_n$
are so large that the estimate cannot be computed in practice.

\subsection{Main results}
\label{se1sub8}
In this article
we define neural network estimates with logistic
squasher activation function
where a linear combination of $K_n$ fully connected
neural networks of depth $L$ and width $r$ is computed in parallel. We use uniform
distributions on the intervals $[-A_n,A_n]$ and $[-B_n,B_n]$ for 
initialization of the input weights and the inner weights,
resp. All outer weights are initially set to zero.
We perform $t_n$ gradient descent steps with stepsize $\lambda_n$,
where we choose
\[
\lambda_n = \frac{1}{\hat{t}_n}
\quad \mbox{and} \quad
t_n = \min \left\{
\hat{t}_n, \lceil (\log n)^{c_8} \cdot K_n^3 \rceil
\right\}
\]
such that
\[
\hat{t}_n \in \left\{
2^i \cdot t_{min} \quad : \quad i \in \N_0 
\right\}
\]
satisfies (with high probability) the following three conditions:
\[
\frac{1}{t_n}
\cdot \sum_{t=0}^{t_n -1}
\lambda_n \cdot \left\|
\nabla_\bw F_n (\bw^{(t)})
\right\|^2
\leq
\frac{c_9}{n},
\]
\[
F_n( \bw^{(t_n)})
\leq
\frac{1}{t_n}
\cdot \sum_{t=0}^{t_n -1}
F_n (\bw^{(t)}) + \frac{c_9}{n},
\]
and
\[
 \max_{t=1, \dots, t_n} \| \bw^{(0)}-\bw^{(t)}\|^2
  \leq
  \frac{c_9 \cdot \log n}{n}.
  \]
  We propose an algorithm which chooses $\lambda_n$ and $t_n$ such that
  these three conditions are satisfied (with high probability). We
  show that the truncated version of the corresponding estimate satisfies
  \begin{eqnarray*}
  &&
  \EXP \int | m_n(x)-m(x)|^2 \PROB_X (dx)
  \\
  &&
\leq
c_{10} \cdot \Bigg(
\EXP \left\{
\inf_{\bw: \|\bw-\bw^{(0)}\| \leq \frac{1}{n}}
\int | f_\bw (x)-m(x)|^2 \PROB_X (dx)
\right\}
+
\frac{A_n^d \cdot B_n^{(L-1) \cdot d}}{n^{1-\epsilon}}
\Bigg)
\end{eqnarray*}
  and, in case of $(p,C)$ smooth regression function  
  and $A_n$ and $B_n$ chosen suitably,
  \[
\EXP \int | m_n(x)-m(x)|^2 \PROB_X (dx)
\leq c_{11} \cdot n^{- \frac{2p}{2p+d} + \epsilon},
  \]
  where $\epsilon>0$ is an arbitrary small number and where $c_{10}, c_{11}>0$
  are constants depending on $\epsilon$.
  Furthermore,
  we implement this estimate and study its finite sample size performance
  in an univariate regression estimation problem, showing that it achieves
  a good performance
  on simulated data. In particular, we observe that for our simulated
  data the above mentioned algorithm for the choice of $\lambda_n$ and
  $t_n$ leads to an estimate which can be computed in a reasonable
  time.

  Our main contributions can be summarized as follows: Motivated
  by a theoretical analysis of the expected $L_2$ error of
  a neural network regression estimated learned by gradient
  descent we propose a special topology of the neural networks,
  a special initialization (where we use uniform distributions
  whose parameters can be considered as smoothing parameters),
  and a special way to choose the stepsize and the number of gradient
  descent steps. We derive a theoretical bound on the expected $L_2$
  error of this estimate, and propose an algorithm where all
  parameters are chosen data-dependent which leads to an estimate
  which outperforms the traditional regression estimates (including
  neural networks estimates) on simulated data in an univariate case. 
  This shows that theoretical analysis of deep neural
  network estimates can lead to new estimates which have an
  improved performance on simulated data.

\subsection{Discussion of related results}
\label{se1sub9}

The huge success of deep learning in applications
has motivated many researchers to investigate theoretically
why these methods are so successful. This has been
studied, e.g., in approximation theory, where quite a few
results concerning the approximation of smooth functions
by deep neural networks have been derived, see
Yarotsky (2018), Yarotsky and Zhevnerchute (2019),
Lu et al. (2020), Langer (2021) and the literature cited therein.
Here it is investigated what kind of topology
and how many nonzero weights are necessary to approximate
a smooth function up to some given error. In applications,
the functions which one wants to approximate has to be estimated
from observed data, which usually contain some random error.
One interesting question in this context is
how well a neural network learned
from such noisy data generalizes on a new independent test data.
Classically this is done within the framework of the VC theory,
and here e.g. the result of
Bartlett et al. (2019) can be used
to bound the VC dimension of classes
of neural networks. For over-parametrized deep
neural networks
(where the number of free parameters adjusted
to the observed data set is much larger than the sample size)
the analysis of the generalization error can be done by using
bounds on the Rademacher complexity
(cf., e.g., Liang, Rakhlin and Sridharan (2015), Golowich, Rakhlin and
Shamir (2019), Lin and Zhang (2019),
Wang and Ma (2022)
and the literature cited therein).
By combining these results it was possible to analyze the
error of least squares regression estimates. Here results
have been shown which indicate why deep learning performs well
in high-dimensional applications: they show
that least squares regression estimates based on deep neural networks
can achieve a dimension reduction in case that the function to be
estimated satisfies a hierarchical composition model, i.e., in case
that it is a composition of smooth functions which do either depend
only on a few components or are rather smooth.
One of the first results in this respect was shown in Kohler and Krzy\.zak (2017),
and later extended by Bauer and
Kohler (2019), Schmidt-Hieber (2020) and Kohler and Langer (2021).
The main trick in these papers is the use of the
network structure of deep neural networks, which implies that the composition
of neural networks is itself a deep neural network. Consequently, any
approximation result for some functions by deep neural networks can be extended 
to approximation of a composition of such functions 
by a deep neural network representing a composition of the approximating
networks. Since in this setting neither the number of weights nor
the depth of the network, which determine the VC dimension and hence the
complexity of the neural network in case that it is not over-parametrized
(cf., Bartlett et al. (2019)), changes much, these neural networks
share the approximation properties and the complexity of neural networks
for low dimensional predictors and hence can achieve dimension
reduction. Bhattacharya, Fan and Mukherjee (2025) showed that a suitably
defined least squares neural network estimate can also achieve
(up to logarithmic factors) optimal rate of convergence results
in interaction models with diverging dimensions.

In practice, least squares estimates cannot be applied because
the corresponding optimization problem cannot be 
solved efficiently. Instead, gradient descent is used to compute
the estimate, and then it is natural to investigate theoretically
whether estimates learned by gradient descent have nice properties.
It was shown in a series of papers,
cf., e.g., Zou et al. (2018), Du et al. (2019),
Allen-Zhu, Li and Song (2019) and Kawaguchi and Huang (2019),
that the application of gradient descent to over-parameterized
deep neural networks can lead to neural networks which (globally)
minimize the empirical risk considered. Unfortunately, as was shown
in  Kohler and Krzy\.zak (2021),
the corresponding estimates do not behave well on a new independent
data.

In applications it is essential to control the approximation, generalization and optimization errors 
simultaneously
(cf., Kutyniok (2020)). Unfortunately, none of the results mentioned above
controls all these three aspects simultaneously.

One way to study these three aspects simultaneously is to use
some equivalent model of deep learning.
The most prominent  approach here is the
neural tangent kernel setting, which was proposed
by Jacot, Gabriel and Hongler (2020).
In this approach  a kernel estimate is studied 
 and its error is used
 to bound the error of the neural network estimate (see also
 Hanin and Nica (2019) and the literature
 cited therein).
 It was observed by Nitanda and Suzuki (2021) that
 in most studies in the
 neural tangent kernel setting the equivalence to deep neural networks
 holds only pointwise and not for the global $L_2$ error, which is
 crucial for predictions problems in practice. So from results
 derived in the neural tangent kernel setting it is often not clear
 how the $L_2$ error of the deep neural network
 estimate behaves.
 An exception is the article
 Nitanda and Suzuki (2021), where the global error
of an over-parametrized shallow neural network
learned by gradient descent was studied based on the neural tangent kernel
approach. However, due to the
use of the neural tangent kernel, the smoothness assumption
on the function to be estimated has to be defined with the aid of
a norm involving the kernel, which does not lead to classical
smoothness conditions usually considered, which makes it hard to
interpret the obtained results. In addition, it is 
required that the number  of neurons be sufficiently large,
but it was not specified what this exactly means, i.e.,  it is
not clear whether the number of neurons must grow e.g.
exponentially in the sample size or not.

Another approach where the estimate
is studied  in some asymptotically equivalent model
is the  mean field approach, cf., e.g.,
Mei, Montanari, and Nguyen (2018), Chizat and Bach (2018),
Nguyen and Pham (2020),
Ba et al. (2020),
Arous, Gheissari and Jagannath (2021),
Bietti et al. (2022),
and the literature cited therein.
Here it is again unclear
how close the behaviour of the deep networks in
the equivalent model is to the behaviour of the deep networks in the
applications, because the equivalent model is based on some approximation of the
deep neural networks using, e.g., some asymptotic expansions.

In a online stochastic gradient setting, where in each gradient descent
step a new independent data point is given, 
Abbe, Adsera, and Misiakiewicz (2023) studies the rate of convergence
of a shallow neural network estimate learned by the layerwise gradient
descent for special regression functions. Here upper and lower bounds
on the rate of convergence (or more precisely: the number of gradient
descent steps required to achieve a given error bound) are derived.

The results presented in this paper are based on the 
statistical theory for deep neural networks developed by the authors
together with various co-authors, see, e.g.
Braun et al. (2023), Drews and Kohler (2023, 2024),
Kohler and Krzy\.zak (2022, 2023) and
Kohler (2024). Here Braun et al. (2023) investigates the rate
for convergence of a shallow neural network estimate learned
by gradient descent. All other papers consider deep neural networks
with the same kind of topology used in the current paper. Kohler
and Krzy\.zak (2023) uses over-parametrized deep ReLU neural network
learned by gradient descent.
Due to the use of Rademacher complexity to control the generalization
error the rate of convergence derived in case of a $(p,C)$--smooth
regression function is of the order $n^{-p/(2p+d) + \epsilon}$ instead
of $n^{-2p/(2p+d) + \epsilon}$ as in the current paper. Drews and Kohler (2024)
derives result concerning the consistency of the estimates, and in
Kohler and Krzy\.zak (2022) and in
Drews and Kohler (2023) the same rate as in the current paper is shown
but only for the special case $p=1/2$.
Here Kohler and Krzy\.zak (2022) use an additional regularization
of the estimate, and Drews and Kohler (2023) shows that this regularization
is not necessary. For general $p$ the above rate
of convergence is derived in Kohler (2024) (again without additional
regularization). Our paper is closely based
on the approach there and shows that the rate of convergence there
can be also achieved with an estimate which uses a data-dependent
choice of the number of gradient descent steps which is in applications
much smaller than the number of gradient descent steps required
in the theoretical result in Kohler (2024).

\subsection{Notation}
\label{se1sub10}
  The sets of natural numbers and real numbers
  are denoted by $\N$ and $\R$, respectively.  For $z \in \R$, we denote
the smallest integer greater than or equal to $z$ by
$\lceil z \rceil$.
The Euclidean norm of $x \in \Rd$
is denoted by $\|x\|$. For a closed and convex set $A \subseteq \R^d$
we denote by $Proj_A x$ that element $Proj_A x \in A$ with
\[
\|x-Proj_A x\|= \min_{z \in A} \|x-z\|.
\]

\subsection{Outline}
\label{se1sub9}
The newly proposed deep learning regression estimate is introduced
in Section \ref{se2}. Section \ref{se3} presents theoretical
results
concerning its rate of convergence. Its finite
sample size performance is illustrated in Section \ref{se4}.
Section \ref{se5} contains the proofs.

\section{Definition of the estimate}
\label{se2}
In the sequel we will
use the logistic squasher (\ref{inteq6}) as activation function.

\subsection{Topology of the network}
\label{se2sub1}
We
let $K_n, L, r \in \N$ be parameters of our estimate and using
these parameters we
 set
\begin{equation}\label{se2eq1}
f_\bw(x) = \sum_{j=1}^{K_n} w_{j,1,1}^{(L)} \cdot f_{j,1}^{(L)}(x)
\end{equation}
for some $w_{1,1,1}^{(L)}, \dots, w_{K_n,1,1}^{(L)} \in \mathbb{R}$, where
$f_{j,1}^{(L)}=f_{\bw,j,1}^{(L)}$ are recursively defined by
\begin{equation}
  \label{se2eq2}
f_{k,i}^{(l)}(x) = 
f_{\bw,k,i}^{(l)}(x) = 
\sigma\left(\sum_{j=1}^{r} w_{k,i,j}^{(l-1)}\cdot f_{k,j}^{(l-1)}(x) + w_{k,i,0}^{(l-1)} \right)
\end{equation}
for some $w_{k,i,0}^{(l-1)}, \dots, w_{k,i, r}^{(l-1)} \in \mathbb{R}$
$(l=2, \dots, L)$
and
\begin{equation}
  \label{se2eq3}
f_{k,i}^{(1)}(x) = 
f_{\bw,k,i}^{(1)}(x) = 
\sigma \left(\sum_{j=1}^d w_{k,i,j}^{(0)}\cdot x^{(j)} + w_{k,i,0}^{(0)} \right)
\end{equation}
for some $w_{k,i,0}^{(0)}, \dots, w_{k,i,d}^{(0)} \in \mathbb{R}$.

This means that we consider neural networks which consist of $K_n$ fully
connected
neural networks of depth $L$ and width $r$ computed in parallel and compute
 a linear combination of the outputs of these $K_n$ neural
networks.
The weights in the $k$-th such network are denoted by
$(w_{k,i,j}^{(l)})_{i,j,l}$, where
$w_{k,i,j}^{(l)}$ is the weight between neuron $j$ in layer
$l$ and neuron $i$ in layer $l+1$.

\subsection{Initialization of the weights}
\label{se2sub2}
We initialize the weights $\bw^{(0)}=((\bw^{(0)})_{k,i,j}^{(l))})_{k,i,j,l}$ as
follows: We set
\[
(\bw^{(0)})_{k,1,1}^{(L)}=0
\quad (k=1, \dots, K_n),
\]
 we choose $(\bw^{(0)})_{k,i,j}^{(l)}$ uniformly distributed on
$[-B, B]$ if $l \in \{1, \dots, L-1\}$, and we
choose
$(\bw^{(0)})_{k,i,j}^{(0)}$ uniformly distributed on
$[-A, A]$, where $A,B \geq 0$ are parameters of the estimate.
Here the random values are defined such that all components
of $\bw^{(0)}$ are independent.

\subsection{Gradient descent}
\label{se2sub3}
Our aim is to choose the weight vector $\bw$ by minimizing
the empirical $L_2$ risk
\begin{equation}
\label{se2eq4}
F_n(\bw) = \frac{1}{n} \sum_{i=1}^n | f_\bw(X_i)-Y_i|^2
\end{equation}
of $f_\bw$
with respect to $\bw$. 

We do this by using gradient descent: Given the random
starting vector $\bw^{(0)}$ for the weights from Subsection \ref{se2sub2}
we compute
$t_n \in \N$ gradient descent steps
\begin{equation}
\label{se2eq5}
\bw^{(t)}
=
\bw^{(t-1)}
-
\lambda_n \cdot \nabla_{\bw} F_n(\bw^{(t-1)})
\quad
(t=1, \dots, t_n)
\end{equation}
with stepsize $\lambda_n>0$.

\subsection{Choice of the stepsize and the number of gradient descent steps}
\label{se2sub4}
We choose
\[
\lambda_n = \frac{1}{\hat{t}_n}
\quad \mbox{and} \quad
t_n = \min \left\{
\hat{t}_n, \lceil (\log n)^{c_{8}} \cdot K_n^3 \rceil
\right\}
\]
such that
\[
\hat{t}_n \in \left\{
2^i \cdot t_{min} \quad : \quad i \in \N_0 
\right\}
\]
satisfies either
 the following three conditions
\begin{equation}
\label{se2eq6}
\frac{1}{t_n}
\cdot \sum_{t=0}^{t_n -1}
\lambda_n \cdot \left\|
\nabla_\bw F_n (\bw^{(t)})
\right\|^2
\leq
\frac{c_9}{n},
\end{equation}
\begin{equation}
\label{se2eq7}
F_n( \bw^{(t_n)})
\leq
\frac{1}{t_n}
\cdot \sum_{t=0}^{t_n -1}
F_n (\bw^{(t)}) + \frac{c_9}{n}
\end{equation}
and
\begin{equation}
  \label{se2eq8}
\max_{t=1, \dots, t_n} \| \bw^{(0)}-\bw^{(t)}\|^2
  \leq
  \frac{c_9 \cdot \log n}{n}
\end{equation}
simultaneously, or such that
\begin{equation}
  \label{se2eq9}
n \cdot (\log n)^{c_{8}} \cdot K_n^3 \leq
  \hat{t}_n \leq 2 \cdot n \cdot (\log n)^{c_{8}} \cdot K_n^3
\end{equation}
holds. We do this by using Algorithm \ref{alg1} below.

\begin{algorithm}
  \caption{Pseudo code for the choice of the stepsize and the number of gradient descent steps. \label{alg1}}
  \KwData{$(x_1,y_1)$, \dots, $(x_n,y_n)$ \\
    $K$, $L$, $r$, $A$, $B$ \\
    $t_{min}=50$, $t_{max,1}= (\log n)^{c_{8}} \cdot K^3$,
    $t_{max,2}=n  \cdot t_{max,1}$, $c_9=10$ 
  }
    \Begin{
        i=0\\
        \Repeat{
\Big($\frac{1}{t}
\cdot \sum_{s=0}^{t -1}
\lambda \cdot \left\|
\nabla_\bw F_n (\bw^{(s)})
\right\|^2
\leq
\frac{c_9}{n} $
and $
F_n( \bw^{(t)})
\leq
\frac{1}{t}
\cdot \sum_{s=0}^{t -1}
F_n (\bw^{(s)}) + \frac{c_9}{n}
$ and $\max_{s=1, \dots, t} \| \bw^{(0)}-\bw^{(s)}\|^2
  \leq
  \frac{c_9 \cdot \log n}{n}
  $
  \Big) or $t \geq t_{max,2}$
        }{
          $\lambda=\frac{1}{2^i \cdot t_{min}}$ \\
          $t=0$ \\
          $\bw^{(0)}=InitializeWeights(K,L,r,A,B)$\\
          \Repeat{$t \geq min(2^i \cdot t_{min}, t_{max,1})$
            or
            $\frac{1}{2^i \cdot t_{min}}
\cdot \sum_{t=0}^{t -1}
\lambda \cdot \left\|
\nabla_\bw F_n (\bw^{(t)})
\right\|^2
>
\frac{c_9}{n} $
or
$\| \bw^{(0)}-\bw^{(t)}\| >
\frac{c_9}{n} $
          }{
$\bw^{(t+1)}
=
\bw^{(t)}
-
\lambda \cdot \nabla_{\bw} F_n(\bw^{(t)})$

$t=t+1$}
                $i=i+1$
        }
        }
  \KwResult{$f_{\bw^{(t)}}$}
\end{algorithm}

    \bigskip
    
    In Algorithm \ref{alg1} the second and the third condition in the inner
    repeat-until loop imply that the first or the third condition in the
    outer repeat-until loop cannot be satisfied if we continue the
    inner loop and therefore the inner loop is terminated
    if one of these conditions holds.

\subsection{Definition of the estimate}
\label{se2sub5}
For the theoretical analysis we consider 
a truncated version of the neural
network with weight vector $\bw^{(t_n)}$, i.e., we define the
estimate by
\begin{equation}
\label{se2eq10}
m_n(x)= T_{\beta_n} (f_{\bw^{(t_n)}}(x))
\end{equation}
 where $\beta_n = c_{12} \cdot \log n$ and $T_{\beta} z
= \max\{ \min\{z, \beta\}, - \beta\}$ for $z \in \R$
and $\beta>0$.

\section{Rate of convergence}
\label{se3}

In this section we present our theoretical results concerning
the estimate introduced in Section \ref{se2}.

\subsection{A general result}

Our first result is a general bound on the expected $L_2$ error
of our estimate.

\begin{theorem}
  \label{th1}
Let $n \in \N$,
let $(X,Y)$, $(X_1,Y_1)$, \dots, $(X_n,Y_n)$
be independent and identically distributed $\Rd \times \R$--valued random variables such that $supp(X)$ is bounded,
the regression function is bounded in absolute value, and
\begin{equation}
\label{th1eq1}
\EXP\left\{
e^{c_7 \cdot Y^2}
\right\}
< \infty
\end{equation}
holds. Let $K_n \in \N$ be such that
\begin{equation}
\label{th1eq2}  
  \frac{K_n}{n^\kappa} \rightarrow 0 \quad (n \rightarrow \infty)
\end{equation}
for some $\kappa>0$, set $A=A_n$ and $B=B_n$ for some
\begin{equation}
\label{th1eq3}  
1 \leq A_n \leq n \quad \mbox{and} \quad
  1 \leq B_n \leq c_{13} \cdot \log n,
\end{equation}
set $\beta_n= c_{12} \cdot \log n$
and define the estimate as in Section \ref{se2}.
Assume $c_7 \cdot c_{12} \geq 3$ and $c_8 > 2L$.
Then we have for any $\epsilon>0$
\begin{eqnarray*}
  &&
  \EXP \int | m_n(x)-m(x)|^2 \PROB_X (dx)
  \\
  &&
\leq
c_{14} \cdot \Bigg(
\EXP \left\{
\inf_{\bw: \|\bw-\bw^{(0)}\| \leq \frac{1}{n}}
\int | f_\bw (x)-m(x)|^2 \PROB_X (dx)
\right\}
+
\frac{A_n^d \cdot B_n^{(L-1) \cdot d}}{n^{1-\epsilon}}
\Bigg)
.
\end{eqnarray*}
\end{theorem}

\noindent
    {\bf Remark 1.}
    The right-hand side above is a sum of two terms. The first
    term
    \[
\EXP \left\{
\inf_{\bw: \|\bw-\bw^{(0)}\| \leq \frac{1}{n}}
\int | f_\bw (x)-m(x)|^2 \PROB_X (dx)
\right\}
    \]
    can be considered as the approximation error of the estimate
    and describes how well the unknown regression function can be
    approximated by deep neural networks whose inner weights are close
    to the randomly initialized weights at the beginning of the
    gradient descent. The second term
    \[
\frac{A_n^d \cdot B_n^{(L-1) \cdot d}}{n^{1-\epsilon}}
    \]
    can be considered as the estimation error of the estimate. It is related
    to the fact that we use gradient descent to minimize the empirical
    $L_2$ risk of the estimate (i.e., the empirical $L_2$ risk on the training
    data) and not the $L_2$ risk.
    
\subsection{Rate of convergence in case of a $(p,C)$--smooth regression function}

If we impose some smoothness condition on the regression function
we can derive an upper bound on the approximation error of the
estimate and use it to derive a bound on the rate of convergence
of the estimate. Our main result in this respect is the following
corollary to Theorem \ref{th1}.

\begin{corollary}
  \label{co1}
  Let $n \in \N$,
let $(X,Y)$, $(X_1,Y_1)$, \dots, $(X_n,Y_n)$
be independent and identically distributed $\Rd \times \R$--valued random variables such that $supp(X)$ is bounded and that (\ref{th1eq1})
holds for some $c_7>0$. Let $p,C>0$ where $p=q+\beta$ for some $q \in \N_0$
and $\beta \in (0,1]$,
  and assume that the regression function
$m:\R^d \rightarrow \R$ is $(p,C)$--smooth.

  Set $\beta_n=c_{12} \cdot \log n$ for some $c_{12}>0$ which satisfies
  $c_7 \cdot c_{12} \geq 3$, and assume $c_8 > 2L$. Let
  $K_n \in \N$ be such that (\ref{th1eq2}) holds for some $\kappa>0$
  and such that
  \[
  \frac{K_n}{n^{175 \cdot (2p+d)^4  \cdot \lceil \log_2(p+d) \rceil}}
  \rightarrow \infty \quad (n \rightarrow \infty)
  \]
  holds.
  Set
  \[
A = A_n = c_{15} \cdot n^{\frac{1}{2p+d}} \cdot \log n \quad \mbox{and} \quad
  B= B_n = c_{16} \cdot \log n,
  \]  
\[
L=\lceil \log_2(q+d) \rceil+1 \quad \mbox{and} \quad r=2 \cdot
\lceil
(2p+d)^2
\rceil
\]
and define the estimate as in Section \ref{se2}.

Then we have for any $\epsilon>0$:
\[
\EXP \int | m_n(x)-m(x)|^2 \PROB_X (dx)
\leq c_{17} \cdot n^{- \frac{2p}{2p+d} + \epsilon}.
\]
  \end{corollary}

\noindent
    {\bf Remark 2.} According to Stone (1982) the optimal minimax
    rate of convergence of the expected $L_2$ error
in case of a $(p,C)$--smooth regression function
    is
    \[
n^{- \frac{2p}{2p+d}}
\]
(cf., e.g., Chapter 3 in Gy\"orfi et al. (2002))
so the rate of convergence above is optimal up to the
arbitrarily small
$\epsilon>0$
    in the exponent. It is an open problem whether a corresponding result
    can also be shown with $\epsilon=0$. In our proof this $\epsilon$
    appears due to our use of the metric entropy bounds for
    bounding the complexity of our over-parametrized space of deep neural
    networks.
    
\section{Application to simulated data}
\label{se4}
In this section we investigate how the estimate behaves on the simulated
data.

For the simulated data we use an example from Gy\"orfi et
al. (2002).
Here we have $d=1$ (so the predictor is univariate) and we choose the
distribution of $X$ to be standard normal restricted to $[-1,1]$,
i.e., the distribution of $X$ has a density which is zero outside
of  $[-1,1]$, and which is proportional to the density
of the standard normal distribution on $[-1,1]$. Then we define
\[
m(x)= \begin{cases}
(x+2)^2/2 & \mbox{if } -1 \leq x < -0.5, \\
x/2 + 0.875 & \mbox{if } -0.5 \leq x < 0, \\
5 \cdot (x-0.2)^2 + 1.075 & \mbox{if } 0 \leq x < 0.5, \\
x+0.125 & \mbox{if }  0.5 \leq x \leq 1, 
\end{cases}
\]
\[
\sigma(x)=0.2-0.1 \cdot \cos ( 2 \cdot \pi \cdot x)
\]
and set
\[
Y=m(X) + \sigma(X) \cdot N
\]
where $N$ is a standard normally distributed random variable
independent of $X$.

We implemented our estimate in $R$ using the logistic squasher activation function, 
and the topology of the network as in
(\ref{se2eq1})-(\ref{se2eq3}), i.e., our network is computing
a linear combination of $K_n$ neural networks with depth $L$ and
width $r$. The initialization is done as described in the previous
section with parameters $A$ and $B$, i.e., all outer weights are
initialized by zero and the weights between the hidden layers
and the weights at the input layer are uniformly distributed on the intervals
$[-B,B]$ and $[-A,A]$, respectively. Then we perform $t_n$ gradient descent
steps with stepsize $\lambda_n$ (maybe adapted to the data 
as described in the previous section).
Here we use the
standard formulas for
backpropagation in order to compute the gradient.

\subsection{Do the estimates generalize well despite
  over-parametrization?}
\label{se4sub1}

We compute our estimate with
parameters $K \in \{100, 200, 400, 800, 1600\}$,
                          $L=4$, $r=8$,  $t_n=K/2$, $\lambda=1/t_n$,
                          $A=1000$
and $B=20$ for $25$ data sets of sample size $n=100$ and compute the
median $L_2$ error and its interquartilerange (IQR). Here the deep neural
network has
\begin{eqnarray*}
  && K*(1+(r+1)+(L-2)*r*(r+1)+r*(d+1)) \\
  &&= K*(1+(8+1)+(4-2)*8*(8+1)+8*(1+1)) = K*170 
\end{eqnarray*}
many weights, so it is clearly over-parametrized.
The median values of the $L_2$ errors and the corresponding IQRs
are reported in Table \ref{se5tab1}.

\begin{table}
\begin{tabular}{|l|l|l|}
\hline
Value of $K$ & number of parameters & median $L_2$ error (IQR)  \\
\hline
$100$ & $17,000$ & $0.0010$ $(0.0597)$            \\
$200$ & $34,000$ & $0.0065$  $(0.0018)$            \\
$400$ & $68,000$ & $0.0039$ $(0.0014)$       \\
$800$ & $136,000$ & $0.0032$  $(0.0010)$            \\
$1600$ & $272,000$ & $0.0036$ $(0.0016)$            \\
\hline
\end{tabular}
\caption{Median $L_2$ errors (and IQRs)
  in $25$  simulations with
  $n=100$,                           $L=4$, $r=8$,  $t_n=K/2$, $\lambda=1/t_n$,
  $A=1000$, $B=20$ and
$K \in \{100,200, 400, 800, 1600\}$  
. \label{se5tab1}}
\end{table}

\begin{figure}[h!]
	  \begin{tabular}{cc}
	    \includegraphics[width=6cm]{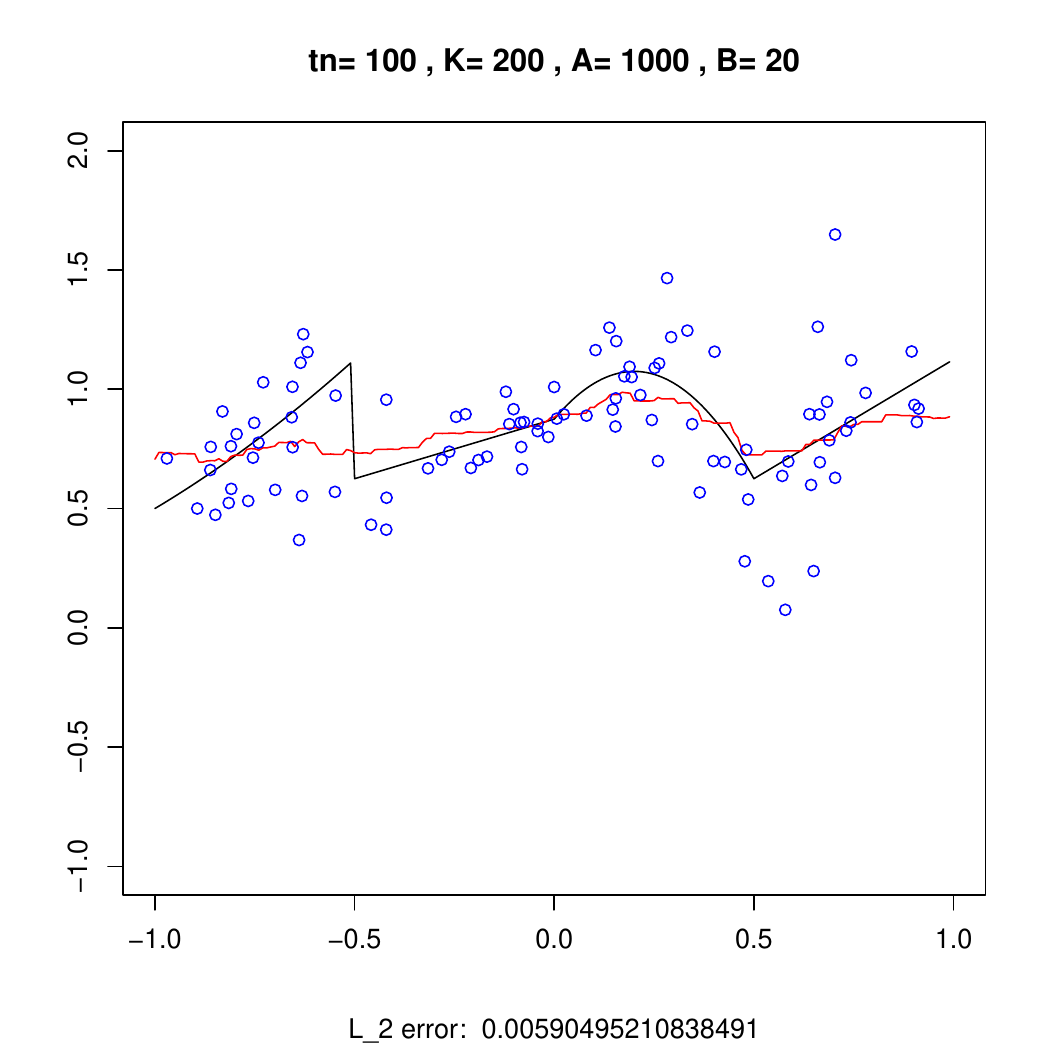}
            &
            		\includegraphics[width=6cm]{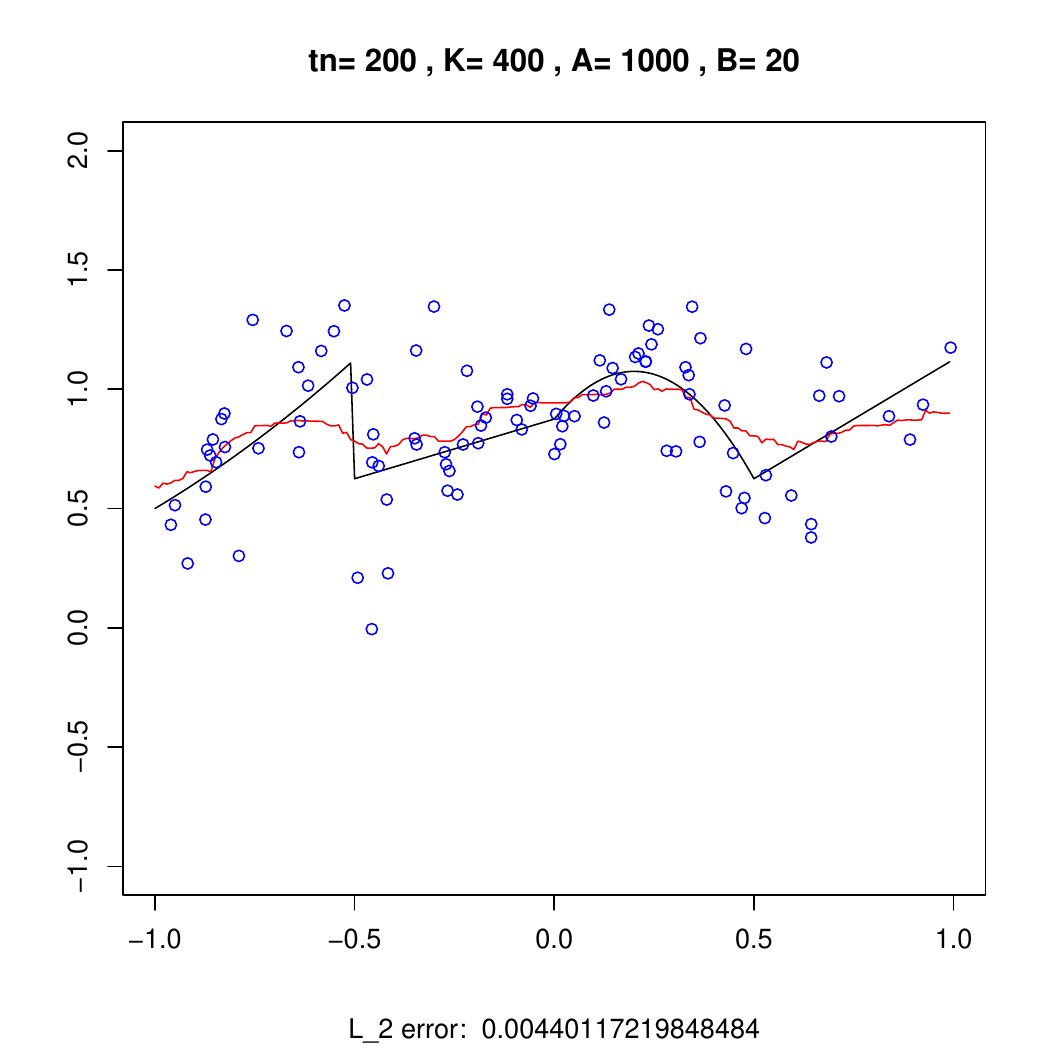}
	                \\
\includegraphics[width=6cm]{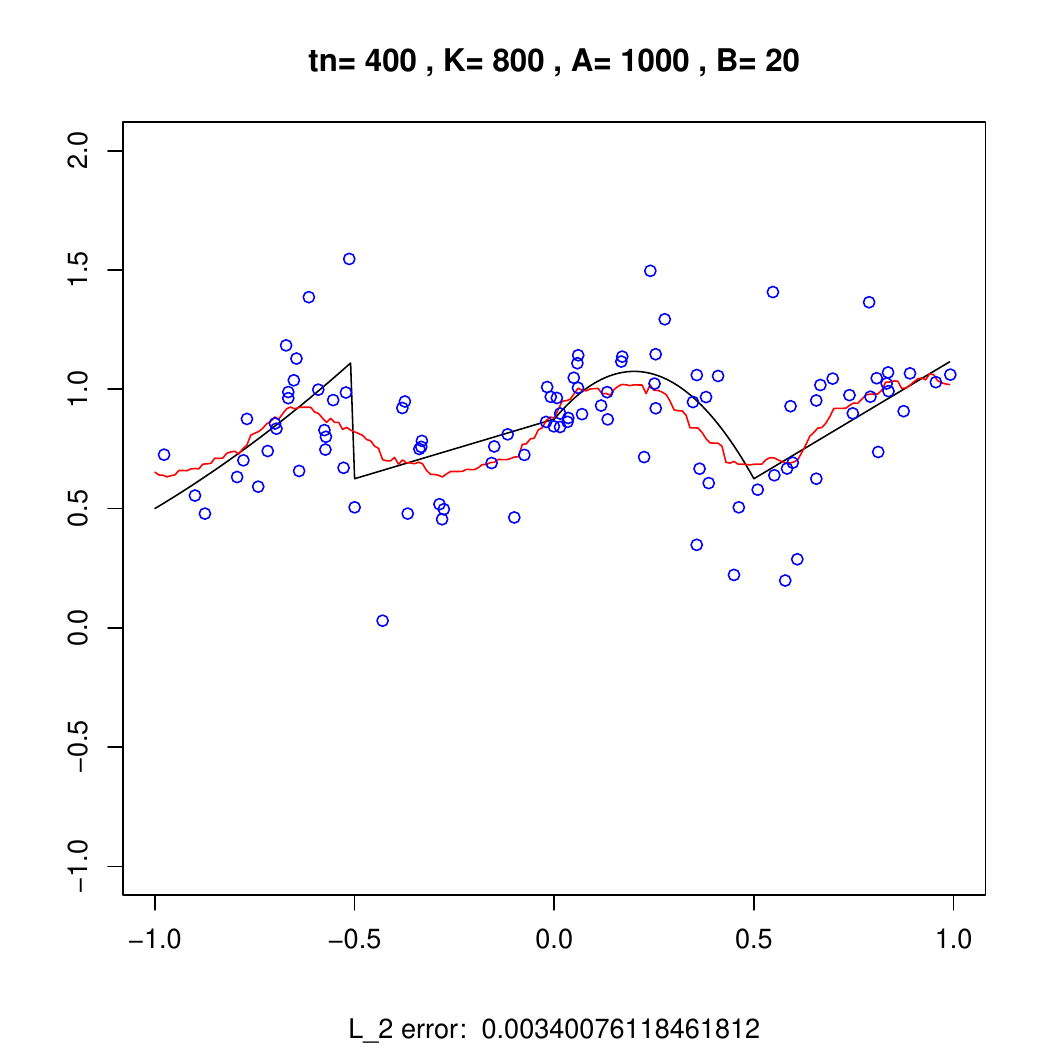}
&
\includegraphics[width=6cm]{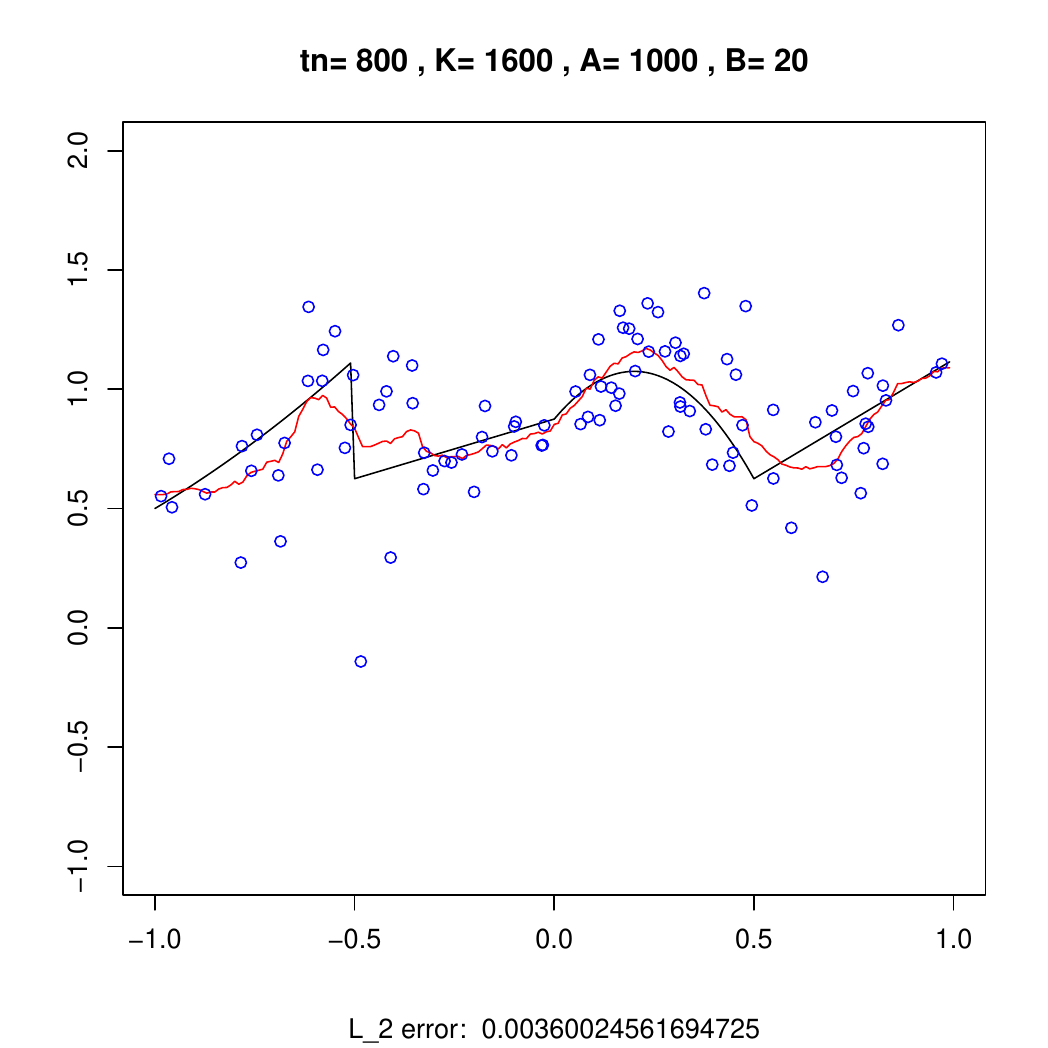}
        \end{tabular}
			\caption{Estimate applied to a sample of size
                          $n=100$, with parameters $K \in \{200, 400, 800, 1600\}$,
                          $L=4$, $r=8$, $\lambda=2/K$, $t_n=K/2$, 
                          $A=1000$ and $B=20$. \label{fig1new}}
	\end{figure}

        Typical estimates for various values of $K$ are shown in Figure \ref{fig1new}.

                As we can see from Table \ref{se5tab1} and Figure \ref{fig1new}
                the error of the estimate decreases with increasing $K$
                as long as $K \leq 800$,
                and although the estimate has much more parameters
                than there are data points, there is no overfitting
                of the data visible even for $K=1600$. We believe that
                the slight increase in the $L_2$ error of the estimate
                for $K=1600$ is either due to the fact that the simple choice
                $t_n=2/K$ is not optimal for a very large value of $K$,
                or that this just occurs because of  random fluctuations of the
                median errors.

\subsection{Are $A$ and $B$ really the smoothing parameters of the estimate?}
\label{se4bsub2}
To see whether the parameters $A$ and $B$ of the uniform distribution are really
the smoothing parameters of the estimate, we apply our estimate with
$n=100$, $K=800$,
                          $L=4$, $r=8$, $\lambda=1/400$, $t_n=400$,
                          $A \in \{ 10, 100, 1000\}$
and $B \in \{2, 20, 200, 2000\}$ to $25$ different data sets and report
the median $L_2$ error and the corresponding IQR in Table \ref{se5tab2}.
\bigskip
\begin{table}
\begin{tabular}{|l|l|l|l|}
\hline
 & $A =10$ & $A =100$ & $A =1,000$   \\
\hline
$B=2$ &  $0.0106$ $(0.0010)$   &  $0.0095$ $(0.0005)$   &  $0.0097$ $(0.0007)$     \\
\hline
$B=20$ &  $0.0034$ $(0.0011)$   &  $0.0033$ $(0.0011)$   &  $0.0034$ $(0.0013)$    \\
\hline
$B=200$ &  $0.0032$ $(0.0018)$   &  $0.0032$ $(0.0010)$   &  $0.0032$ $(0.0022)$       \\
\hline
$B=2000$ &  $0.0035$ $(0.0010)$   &  $0.0030$ $(0.0016)$   &  $0.0034$ $(0.0021)$     \\
\hline
\end{tabular}
\caption{Median $L_2$ errors (and IQRs) in $25$  simulations with
  $n=100$, $K=800$,                          $L=4$, $r=8$,  $t_n=400$, $\lambda=1/400$,
  $A \in \{10, 100, 1000\}$ and $B \in \{2, 20, 200, 2000\}$.
  \label{se5tab2}}
\end{table}

\bigskip

\noindent

We clearly see that $A$ and $B$  have an influence on the $L_2$ error.
If $B$ is too small the $L_2$ errors get large. Otherwise
it is not clear how the values of $A$ and $B$
influence the errors. We think this is due to the fact that
they influence simultaneously the generalization error (where
larger values increase the error) and the approximation error
(where large values of $A$ decrease the approximation error, and where
very
large values of $B$ might decrease the approximation error again
because large values of $A$ might result in input neurons with
an nearly constant output for which a larger value of $B$ might be
an adavantage).

\subsection{Do the data-dependent choices of the stepsize and the number
  of gradient descent steps work?}
\label{se4bsub3}

\begin{table}
  \begin{center}
\begin{tabular}{|c|c|c|}
\hline
Value of $K$  & median $L_2$ error (IQR) & number simulations with $t_n \neq K/2$  \\
\hline
$100$ &  $0.0082$ $(0.0015)$ & $22$            \\
$200$ &  $0.0059$  $(0.0011)$ & $7$           \\
$400$ &  $0.0044$ $(0.0007)$  & $2$     \\
$800$ &  $0.0034$  $(0.0016)$ & $1$           \\
$1600$ &  $0.0034$ $(0.0016)$ & $0$           \\
\hline
\end{tabular}
\end{center}
  \caption{Median $L_2$ errors (and IQRs) in $25$  simulations
    of the adaptive estimate for
    $n=100$, $L=4$, $r=8$, $A=1000$, $B=20$ and
    $K \in \{100, 200, 400, 800, 1600\}$. \label{se5tab3}}
\end{table}

\begin{figure}[h!]
	  \begin{tabular}{cc}
	    \includegraphics[width=7cm]{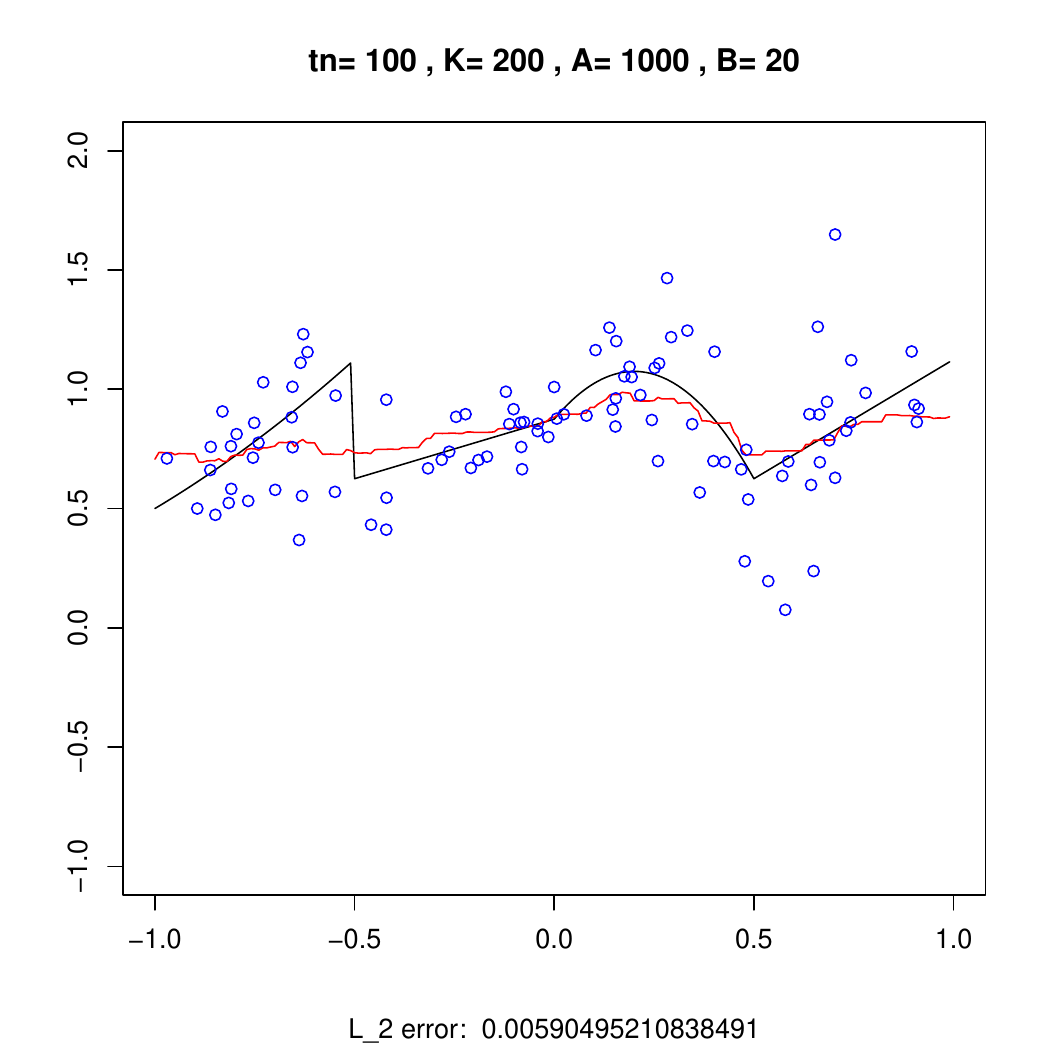}
            &
            		\includegraphics[width=7cm]{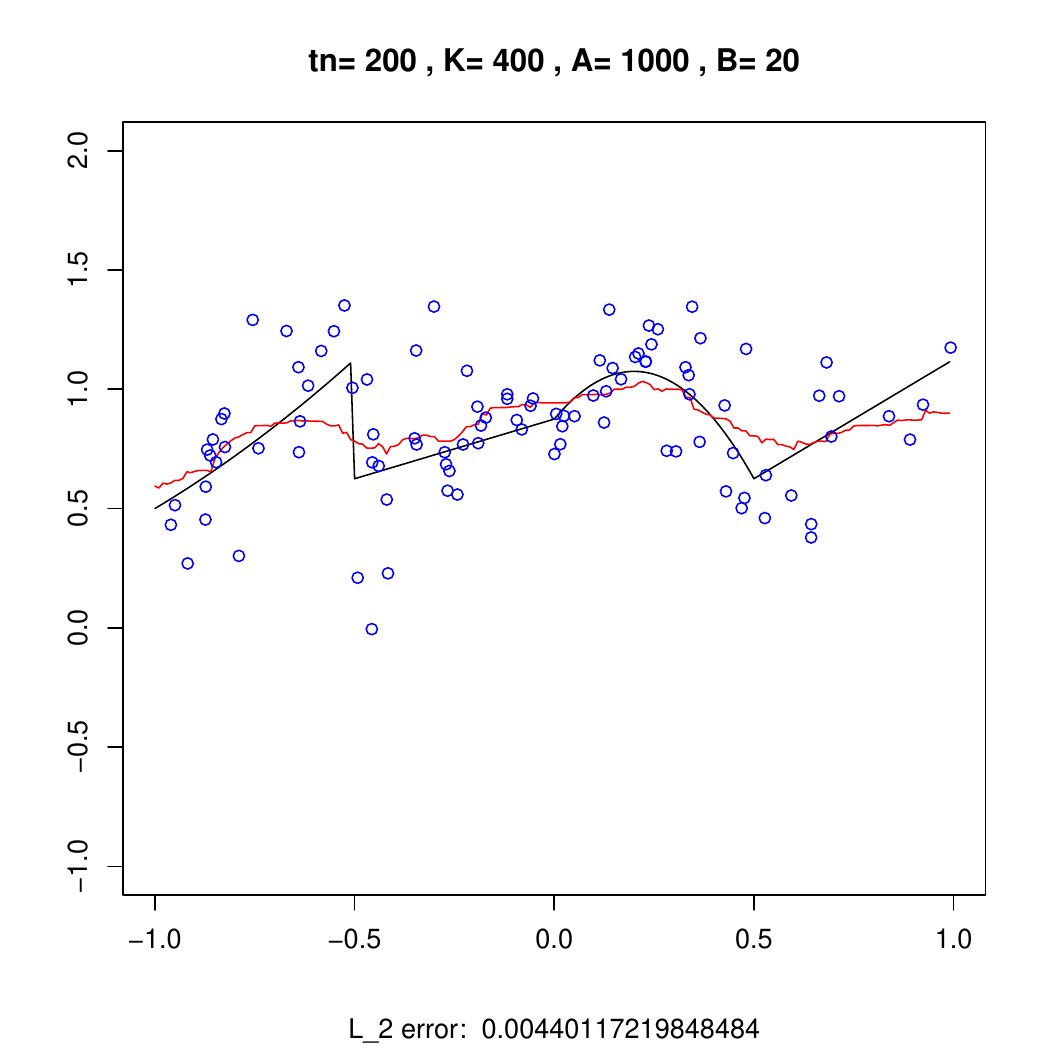}
                        \\
                        	    \includegraphics[width=7cm]{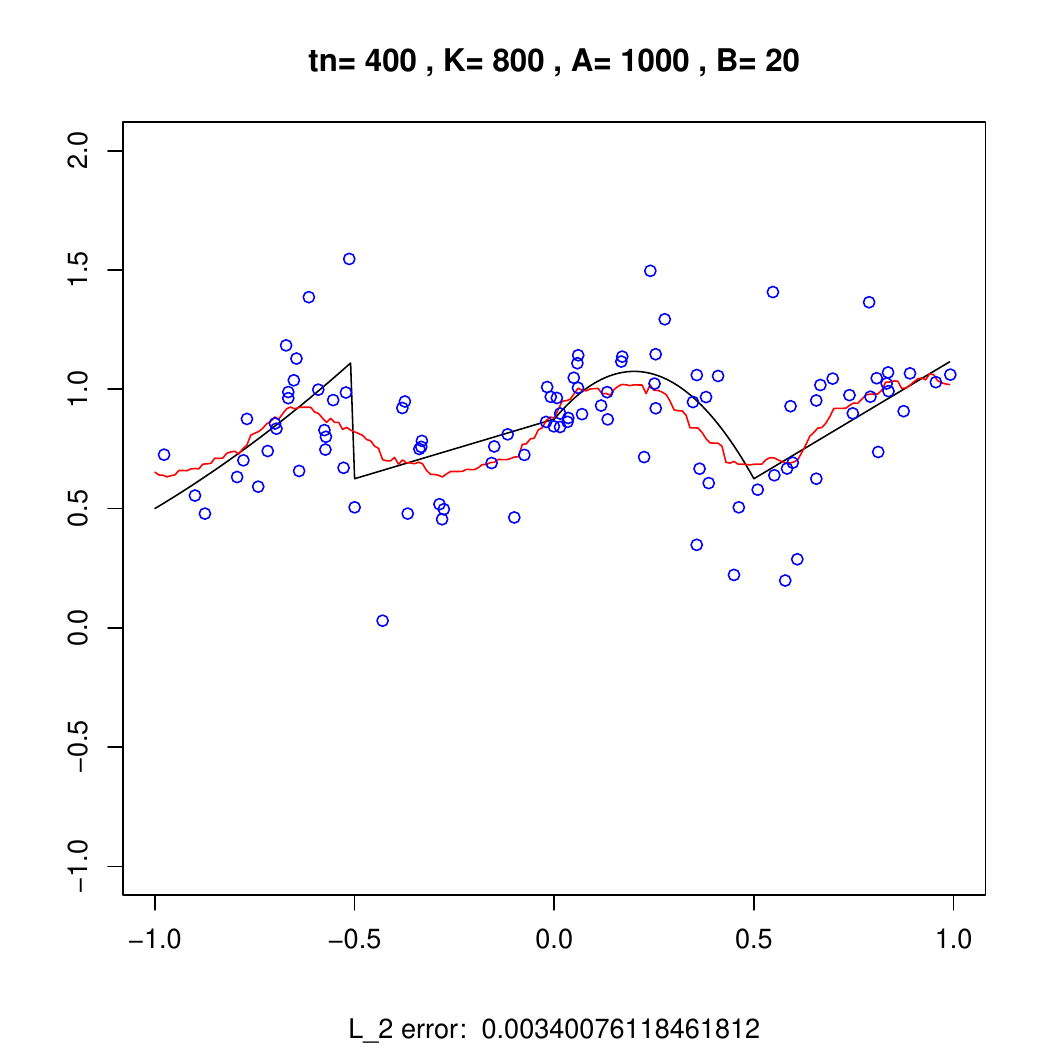}
            &
            		\includegraphics[width=7cm]{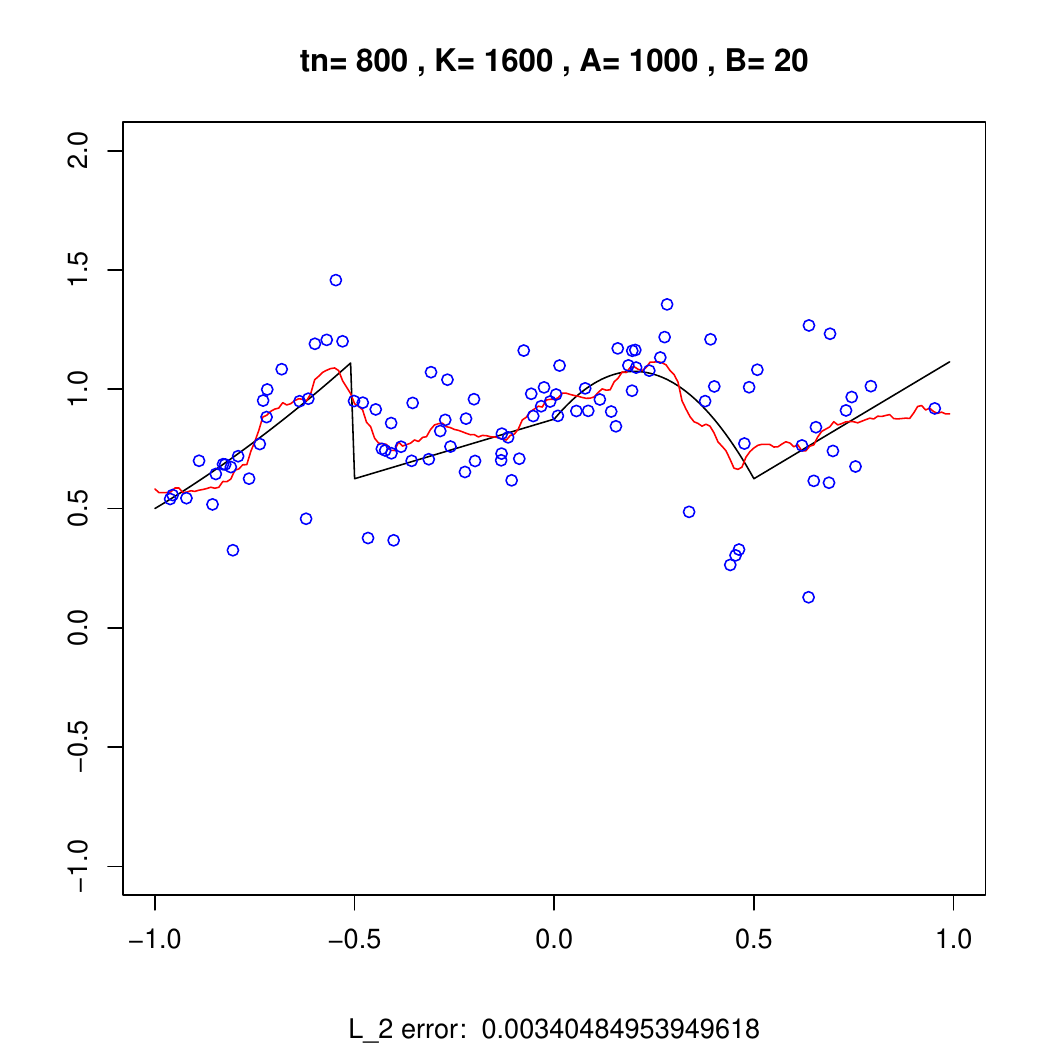}
		                        \\
            		
          \end{tabular}
			\caption{Adaptive estimates applied to a sample of size
                          $n=100$, with parameters $K \in \{200, 400, 800,1600\}$,
                          $L=4$, $r=8$,  
                          $A=1000$ and $B=20$.   \label{fig2}}
	\end{figure}

In this subsection we investigate whether the proposed
data-dependent choice of the stepsize and the number
  of gradient descent steps improves the estimate. To do this, 
  we apply our adaptive estimate, where the number of gradient
  descent steps and the stepsize is chosen as in Subsection
\ref{se2sub4}
  with
$n=100$, $K \in \{100, 200, 400, 800\}$,
                          $L=4$, $r=8$,
                          $A=1000$
  and $B =20$ to $25$ different data sets.
  The median values of
  the $L_2$ errors and their IQRs are reported in Table \ref{se5tab3}.
  There we also report in how many of the $25$ simulations the
  adaptive estimate chooses $t_n \neq K/2$ (and hence uses a different
  value than the non-adaptive estimate in Table \ref{se5tab1}).
   In Figure \ref{fig2}
   we show plots of typical estimates which we get for different
   values of $K$.

                Again the error of the estimate gets smaller with increasing $K$.
                For large values of $K$ the adaptive algorithm always chooses $t_n=2/K$
                which explains why there is no improvement in comparison
                with the non adaptive choice of $t_n$ and $\lambda_n$.
                However, for small values of $K$
                the median errors with the adaptive choice of $\lambda$ and $t_n$ are smaller than the median errors for $t_n=K/2$ and $\lambda_n=1/t_n$
                in Table \ref{se5tab1}.

\subsection{Is an adaptive choice of the weights bounds $A$ and $B$ during
  initialization useful?}
\label{se4bsub4}
We have identified the parameters $A$ and $B$ of the uniform
distributions for the initialization of the weights as possible
smoothing parameters of our neural network estimate.
In this subsection we investigate whether it is useful
to choose these parameters in data-dependent way using
splitting of the sample (cf., e.g., Chapter 7 in Gy\"orfi et
al. (2002)).
Here the given data is divided into the training data consisting of the
first $n_l$ data points, and the testing data consisting of the
$n_t=n-n_l$ remaining data points (e.g., with $n_l \approx n/2$
or $n_l \approx \frac{2}{3} \cdot n$). Then a finite set $\P$ of possible values
for $(A,B)$ is selected, for each value of $(A,B)$ of this set the estimate
\[
m_{n,(A,B)}(\cdot)=m_{n_l,(A,B)}(\cdot, \D_{n_l})
\]
is computed using this value of $(A,B)$ and only the training data, 
and finally that value $(\hat{A},\hat{B}) \in \P$ is selected for which
the empirical $L_2$ risk on the testing data is minimal. Thus, we
compute
\[
(\hat{A},\hat{B}) = \arg \min_{(A,B) \in \P}
\frac{1}{n_t} \sum_{i=n_l+1}^n | Y_i - m_{n_l,(A,B)}(X_i)|^2
\]
and use as estimate
\[
m_n(\cdot) = m_{n_l, (\hat{A},\hat{B})}(\cdot, \D_{n_l}).
\]

We compute this estimate $25$-times with $n=100$, $n_{train}=80$,
$n_{test}=20$, $K \in \{100,200,400,800,1600\}$, $L=4$, $r=8$
and choose the stepsize and the number of gradient descent steps
adaptively as in the previous section and choose the parameters
$A$ and $B$ adaptively from the sets $A \in \{10, 100, 1000\}$
and $B \in \{20, 200, 2000\}$ via splitting of the sample.
The results are reported in Table \ref{se5tab4}.

\begin{table}
  \begin{center}
\begin{tabular}{|l|l|}
\hline
Value of $K$  & median $L_2$ error (IQR)  \\
\hline
$100$ &  $0.0069$ $(0.0016)$            \\
$200$ &  $0.0051$  $(0.0011)$            \\
$400$ &  $0.0046$ $(0.0015)$       \\
$800$ &  $0.0036$  $(0.0014)$             \\
\hline
\end{tabular}
\end{center}
  \caption{Median $L_2$ errors (and IQRs) in $25$  simulations
    of the adaptive estimate for
    $n=100$, $L=4$, $r=8$ and
    $K \in \{100, 200, 400, 800\}$, where we choose 
$A \in \{10, 100, 1000\}$
and $B \in \{20, 200, 2000\}$ via splitting of the sample. \label{se5tab4}}
\end{table}

\begin{figure}[h!]
	  \begin{tabular}{cc}
    \includegraphics[width=6cm]{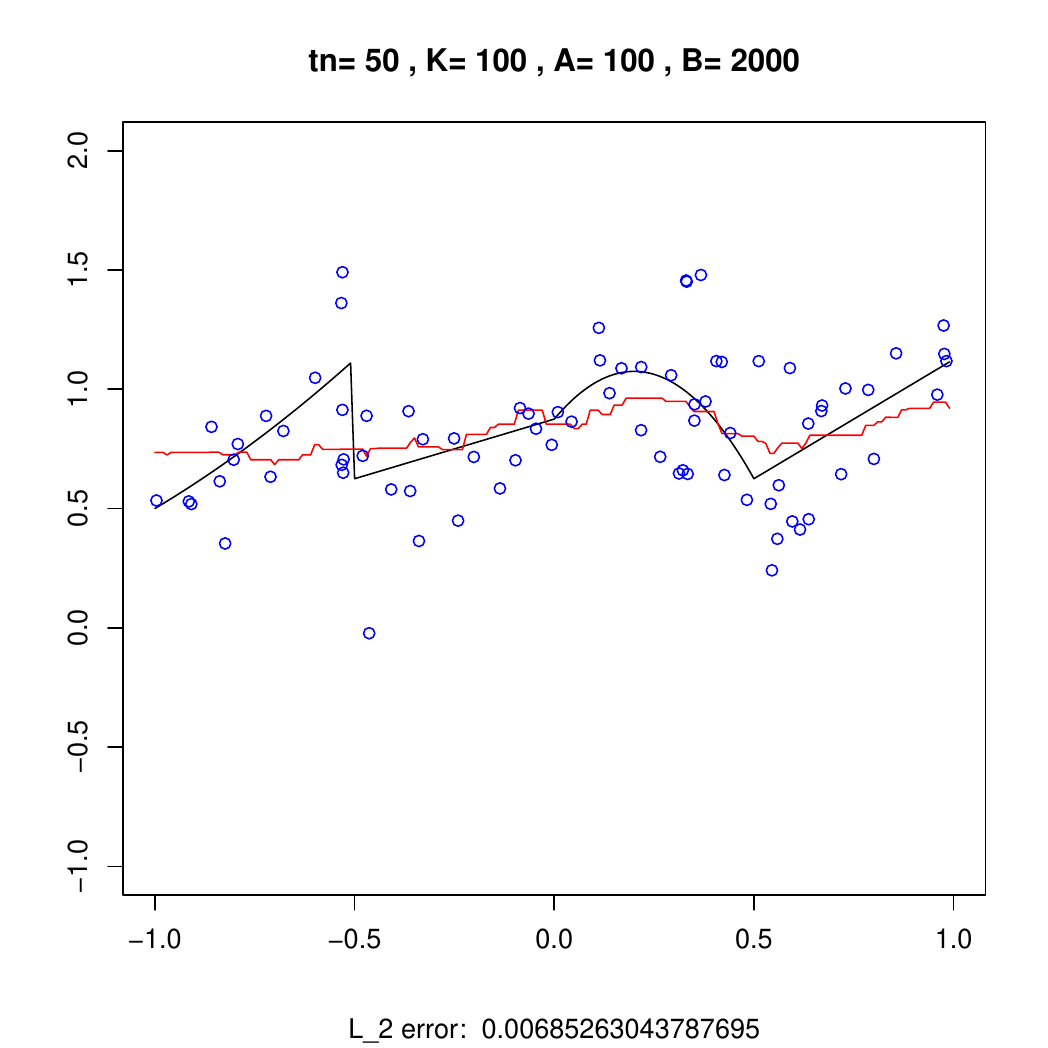}
            &
            		\includegraphics[width=6cm]{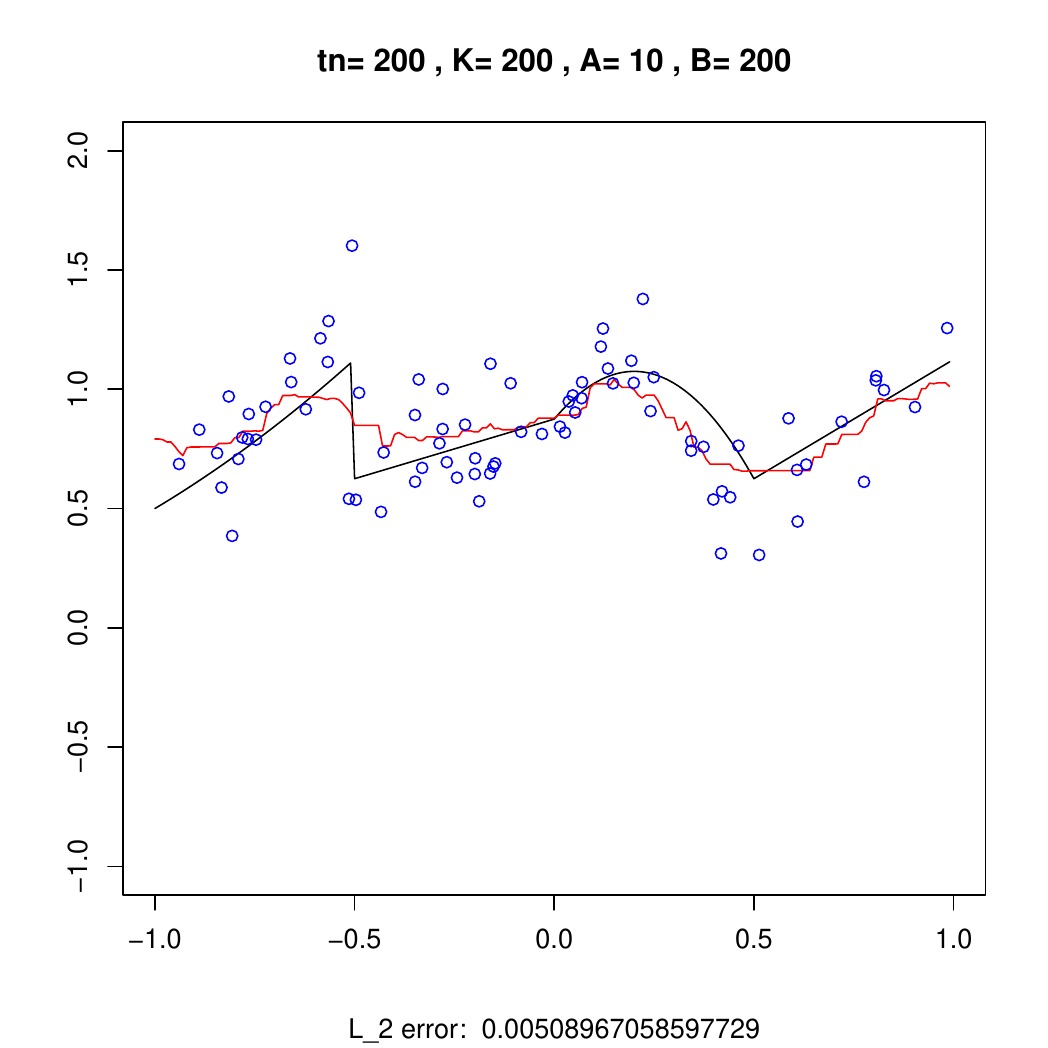}
	                \\
\includegraphics[width=6cm]{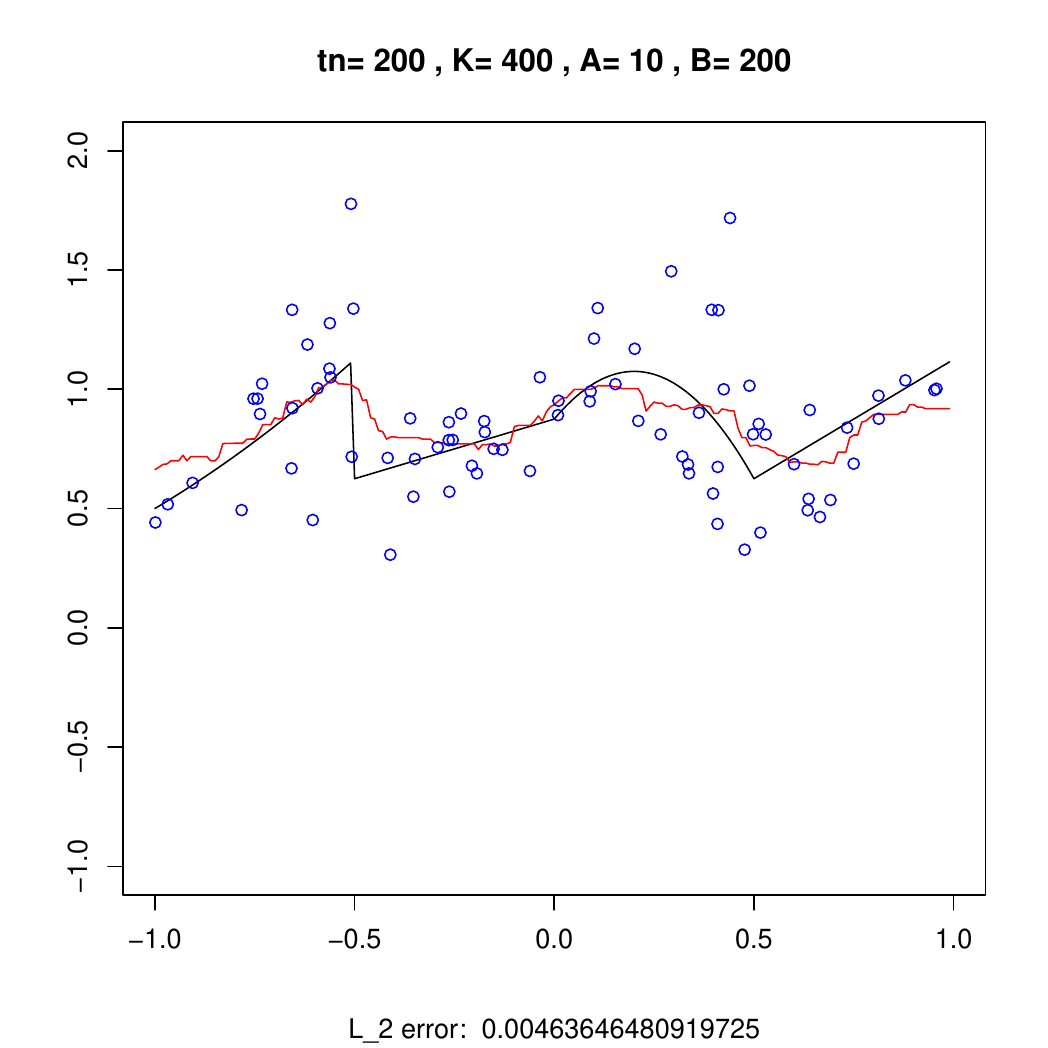}
&
\includegraphics[width=6cm]{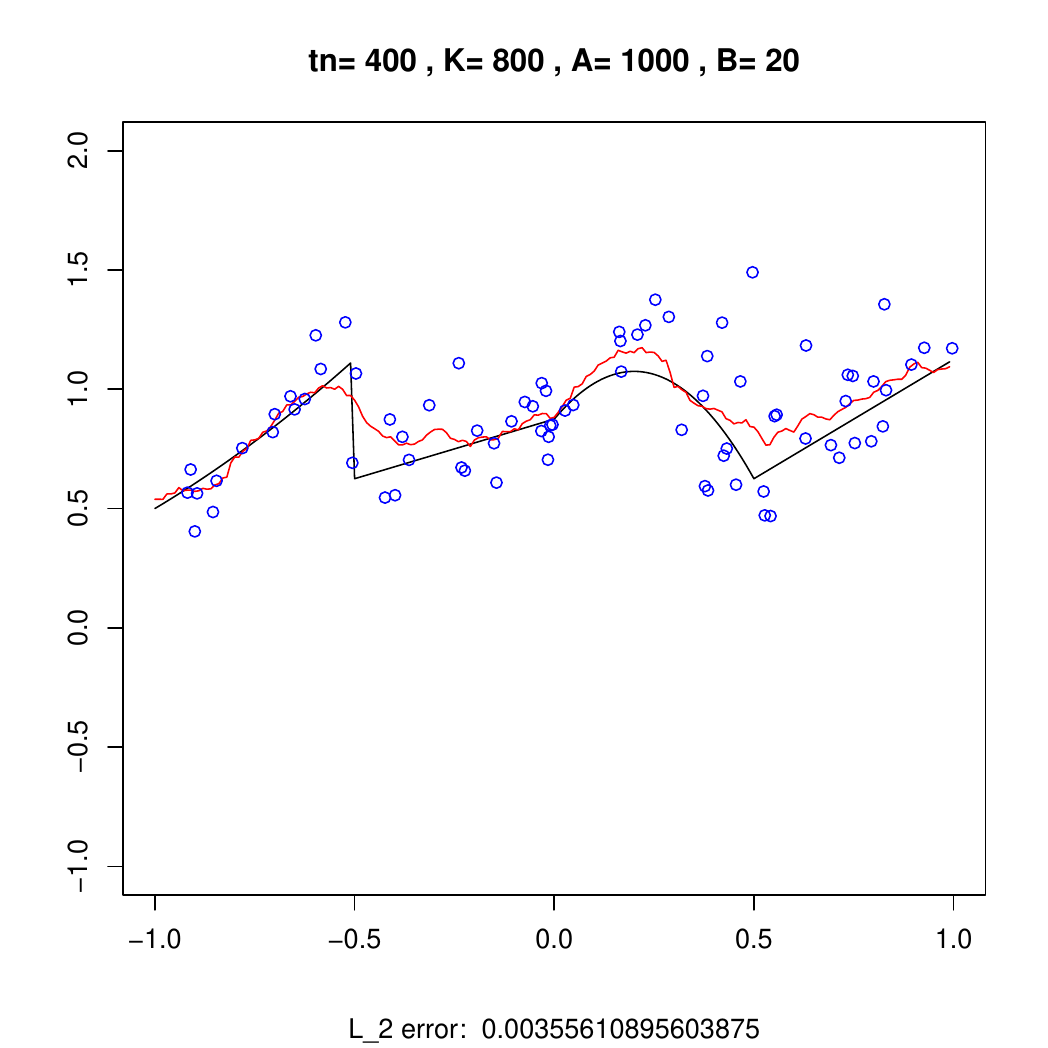}
        \end{tabular}
			\caption{Estimate applied to a sample of size
                          $n=100$, with parameters $K \in \{100, 200, 400, 800\}$,
                          $L=4$, $r=8$, adaptively chosen values for $\lambda$
                          and $t_n$, and values of
$A \in \{10, 100, 1000\}$
                          and $B \in \{20, 200, 2000\}$
                          chosen via splitting of the sample with $n_{train}=80$
                          and $n_{test}=20$.  \label{fig5}}
	\end{figure}

   Plots of typical estimates which we get for the different
   values of $K$ are shown in Figure \ref{fig5}.

   The results show that for small values of $K$ this adaptive estimate
   yields a smaller error than the non-adaptive estimate.
   For large values of $K$ the error of the estimate
   is approximately the same as for the
   other estimates, although
   it is based mainly on only $80\%$ of the data.

\subsection{How good is the estimate?}
\label{se4bsub5}
\begin{figure}[h!]
                  \begin{center}
	  \begin{tabular}{cc}
    \includegraphics[width=6cm]{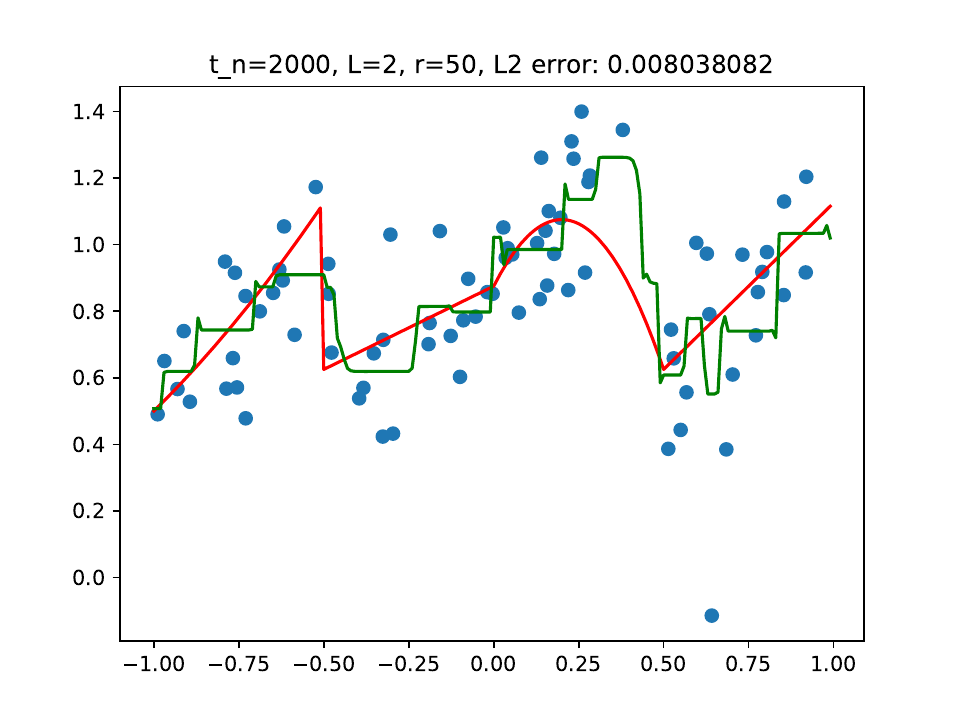}
            &
            		\includegraphics[width=6cm]{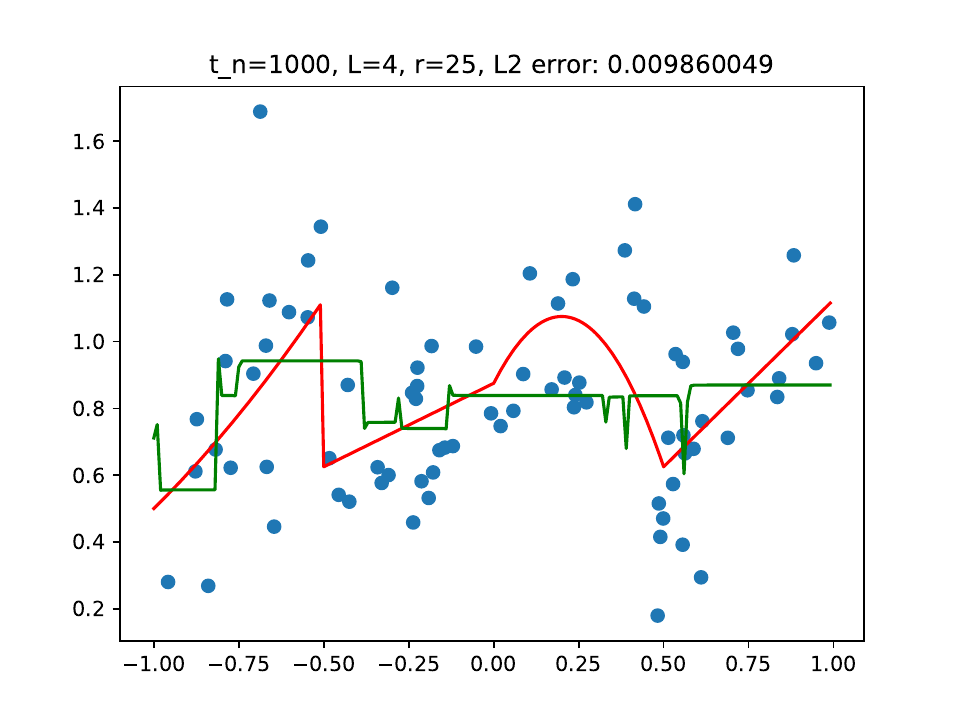}
	                \\
\includegraphics[width=6cm]{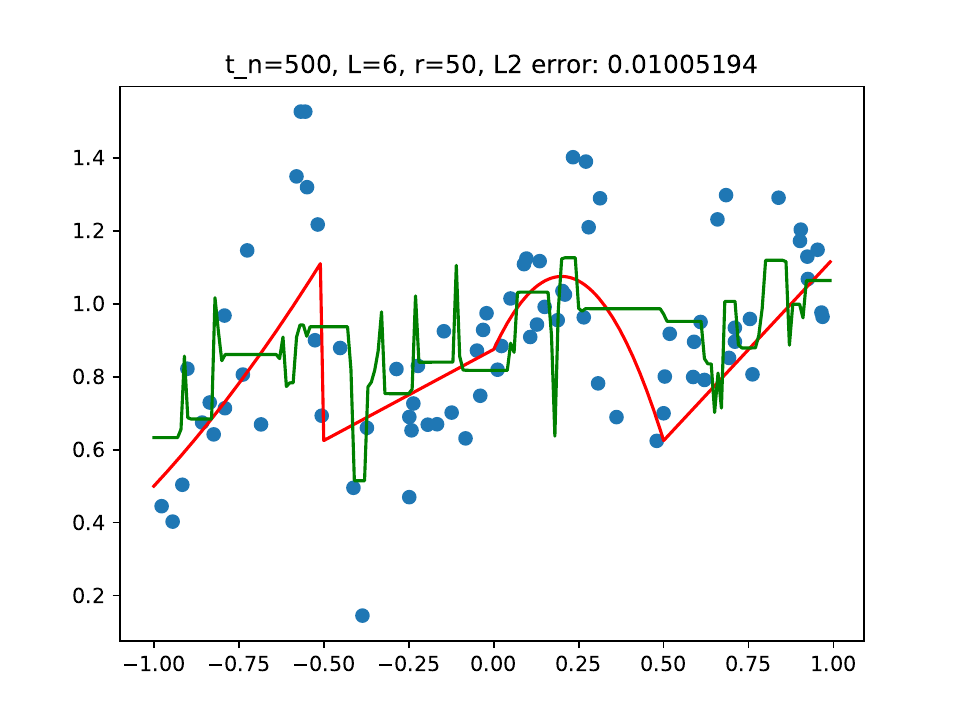}
&
\includegraphics[width=6cm]{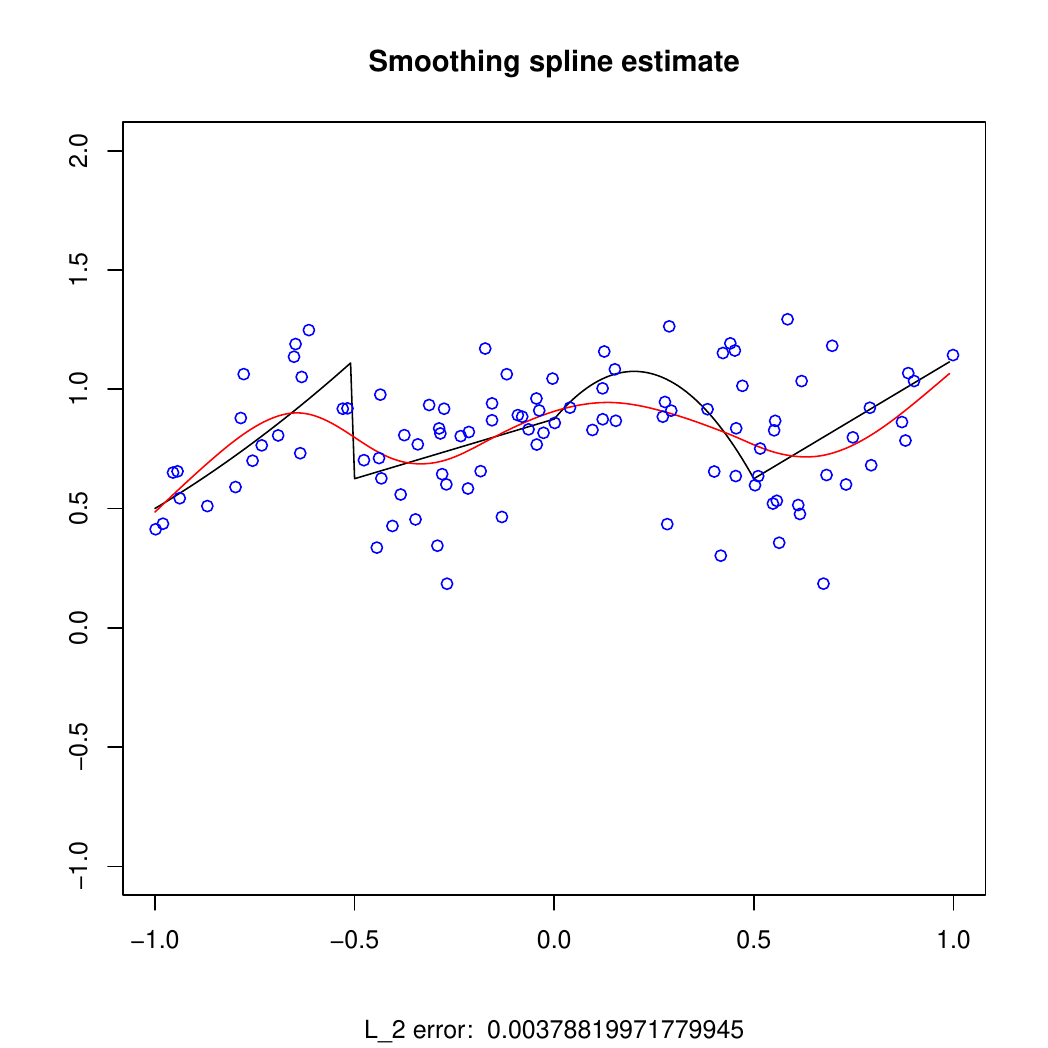}
\\
          \end{tabular}
          \end{center}
			\caption{Standard neural network estimates with $L=2$, $L=4$ and $L=6$ hidden layers and a smoothing spline estimate applied each time to a sample of size
                          $n=100$.  \label{fig6}}
	\end{figure}

In order to see how good our newly introduced neural network regression
estimate
is compared with other known estimates, we apply to our data also
standard neural network estimates with $2$, $4$ and $6$ hidden layers
and a data-dependent chosen number of hidden neurons, and
a smoothing spline estimate. For the neural network estimates
$nnfc2$, $nnfc4$ and $nnfc6$
with $2$, $4$ and $6$ hidden layers, resp.,
the number $r \in \{10, 25, 50, 100, 200\}$ of hidden neurons
and the number $t_n \in \{500, 1000, 2000\}$ of gradient descent steps is chosen data-dependent using splitting of the sample with $n_{train}=80$
and $n_{test}=20$.
The estimate uses the logistic squasher as activation function and the
initialization of the weights is done as before, i.e., all outer weights are
initialized by zero and the weights between the hidden layers
and the weights at the input layer are uniformly distributed on the intervals
$[-20,20]$ and $[-1000,1000]$, respectively. 
The estimates are implemented in Python using the package
{\it tensorflow} with gradient descent as implemented in this package
using the ADAM rule for the data-dependent choice of the stepsize.
The smoothing spline estimate $smooth-spline$ is applied
as as implemented in R by the procedure
{\it Tps()} from the library {\it fields}. The smoothing parameter
of this estimate is chosen data dependent by generalized cross validation
as implemented in {\it Tps()}.
We apply each of these estimates $25$ times to independent data sets of sample size $n=100$. The results are reported in Table \ref{se5tab5}.
\begin{table}
  \begin{center}
\begin{tabular}{|l|l|}
\hline
Estimate  & median $L_2$ error (IQR)  \\
\hline
$nnfc2$ &  $0.0080$ $(0.0030)$            \\
$nnfc4$ &  $0.0099$  $(0.0047)$            \\
$nnfc6$ &  $0.0100$ $(0.0059)$       \\
$smooth-spline$ &  $0.0038$  $(0.0026)$            \\
\hline
\end{tabular}
\end{center}
  \caption{Median $L_2$ errors in $25$  simulations
    of the three different standard neural network estimates
    and the smoothing spline estimate.
 \label{se5tab5}}
\end{table}
   Plots of typical estimates which we get for the different
   estimates are shown in Figure \ref{fig6}.

                We see that the median $L_2$ errors of the neural network estimates
                in Table \ref{se5tab5} are substantially larger than the
                median $L_2$ error of the smoothing spline estimate.
                In contrast, the newly proposed deep neural network estimates
                of this paper achieve for $K \geq 800$
                a performance which is as good or
                even slightly better than this smoothing spline estimate.
                
                This shows that our theoretical
                approach to deep learning improves in our example deep
                neural network estimates drastically such that they become
                comparable good as a standard estimate in an univariate
                regression problem. Of course, in this case the standard
                estimate is much easier to compute, however the potential
                of this result is that by modifying deep neural networks
                in the multivariate case in the same way (which requires
                an extension of the currently available theory for deep
                neural network estimates learned by gradient descent)
                might lead to an improvement of the deep neural network
                estimates in a case where standard estimates do not outperform
                them (because their results in high-dimensional settings
                are not as good as standard
                deep neural network estimates as is shown, e.g.,
                in the simulations in
                Bauer and Kohler (2019)).

\section{Proofs}
\label{se5}

\subsection{An auxiliary result for the proof of Theorem \ref{th1}}

In the proof of Theorem \ref{th1} we will use the following lemma
in order to analyze the gradient descent.

\begin{lemma}
  \label{le1}
  Let $d, J_n \in \N$, and for $\bw \in \R^{J_n}$ let
$
  f_{\bw}: \R^d \rightarrow \R
  $
  be a (deep) neural network with weight vector $\bw$.
  Assume that for each $x \in \R^d$
  \[
\bw \mapsto f_\bw(x)
  \]
is a continuously differentiable function on $\R^{J_n}$.
  Let
\[
F_n(\bw) = \frac{1}{n} \sum_{i=1}^n |Y_i - f_{\bw}(X_i)|^2
\]
be the empirical $L_2$ risk of $f_{\bw}$, and use gradient
descent in order to minimize $F_n(\bw)$. To do this, choose
a starting weight vector $\bw^{(0)} \in \R^{J_n}$, choose
$\delta_n \geq 0$ and let
\[
A \subset \left\{
\bw \in \R^{J_n} \, : \, \|\bw-\bw^{(0)}\| \leq \delta_n
\right\}
\]
be a closed and convex set of weight vectors. Choose a stepsize
$\lambda_n > 0$ and a number of gradient descent steps $t_n \in \N$
and compute
\[
\bw^{(t+1)} = Proj_A \left(
\bw^{(t)}
-
\lambda_n
\cdot \nabla_{\bw} F_n( \bw^{(t)} )
\right)
\]
for $t=0, \dots, t_n-1$.

Let $C_n \geq 0$, $\beta_n \geq 1$ and assume
\begin{equation}
  \label{le1eq1}
  \sum_{j=1}^{J_n}
  \left|
  \frac{\partial}{\partial w^{(j)}}
  f_{\bw_1}(x)
-
\frac{\partial}{\partial w^{(j)}}
f_{\bw_2}(x)
\right|^2
   \leq C_n^2 \cdot
  \| \bw_1 - \bw_2 \|^2
  \end{equation}
  for all  $\bw_1, \bw_2 \in A$, $x \in \{X_1, \dots, X_n\}$,
\begin{equation}
  \label{le1eq3}
  |Y_i| \leq \beta_n \quad (i=1, \dots, n) 
\end{equation}
and
\begin{equation}
  \label{le1eq4}
C_n \cdot \delta_n^2 \leq 1.
  \end{equation}
Let $\bw^* \in A$ and assume
\begin{equation}
  \label{le1eq5}
  |f_{\bw^*}(x)| \leq \beta_n \quad (x \in \{X_1, \dots, X_n\}).
\end{equation}
Then  
\begin{eqnarray*}
\frac{1}{t_n}  \sum_{t=0}^{t_n-1} F_n(\bw^{(t)})
  &\leq&
 F_n(\bw^*) + \frac{\| \bw^* - \bw^{(0)} \|^2}{2 \cdot \lambda_n \cdot t_n}
  +
  3 \cdot \beta_n \cdot C_n \cdot  \frac{1}{t_n}  \sum_{t=0}^{t_n-1}
  \| \bw^* - \bw^{(t)} \|^2 \\
  &&
  +
\frac{1}{2} \cdot
  \lambda_n \cdot  
\frac{1}{t_n} \sum_{t=0}^{t_n-1} \|\nabla_\bw F_n(\bw^{(t)}) \|^2.
  \end{eqnarray*}
  \end{lemma}

\noindent
    {\bf Proof.} The result follows by a straightforward
    modification of the proof of Lemma 1 in Kohler (2024).
    For the sake of completeness we nevertheless present
    a more or less complete proof here.
    
    The basic idea of the proof is to analyze the
    gradient descent by relating it to the gradient
    descent of the linear Taylor polynomial of $f_\bw$. To do this,
    we define
    for $\bw_0, \bw \in \R^{J_n}$ the linear Taylor polynomial
    of $f_\bw(x)$ around $\bw_0$ by
    \[
    f_{lin, \bw_0, \bw}(x)
    =
    f_{\bw_0}(x)+
    \sum_{j=1}^{J_n}
    \frac{\partial f_{\bw_0}(x)}{ \partial \bw^{(j)}}
      \cdot
      (\bw^{(j)} - \bw_0^{(j)})
    \]
    and introduce the empirical $L_2$ risk of this
    linear approximation of $f_\bw$ by
    \[
    F_{n,lin,\bw_0}(\bw)
    =
    \frac{1}{n} \sum_{i=1}^n |Y_i - f_{lin,\bw_0,\bw}(X_i)|^2.
    \]
Then
$F_{n,lin,\bw_0}(\bw)$ is as a function of $\bw$ a convex function
(cf. Kohler (2024), proof of Lemma 1).

Because of 
$   f_{lin, \bw_0, \bw_0}(x)
    =
    f_{\bw_0}(x)$ and
    $   \nabla_\bw f_{lin, \bw_0, \bw_0}(x)
    =
    \nabla_w f_{\bw_0}(x)$ we have
    \[
    F_{n,lin,\bw^{(t)}}(\bw^{(t)}) = F_n(\bw^{(t)})
      \quad \mbox{and} \quad
    \nabla_w F_{n,lin,\bw^{(t)}}(\bw^{(t)})
    =
    \nabla_w F_{n}(\bw^{(t)}),
    \]
    hence $\bw^{(t+1)}$ is computed from $\bw^{(t)}$ by one gradient descent
    step
    \[
\bw^{(t+1)} = Proj_A \left(
\bw^{(t)}
-
\lambda_n
\cdot \nabla_{\bw} F_{n,lin,\bw^{(t)}}( \bw^{(t)} )
\right)
    \]
    applied to the convex function $ F_{n,lin,\bw^{(t)}}( \bw )$.
    This will enable us to use techniques for the analysis of the gradient descent
    for convex functions in order to analyze the gradient descent applied
    to the nonconvex function $F_n(\bw)$.

    In order to do this we observe
    \begin{eqnarray*}
      &&
\frac{1}{t_n}  \sum_{t=0}^{t_n-1} F_n(\bw^{(t)})
  -
  F_n(\bw^*)
  \\
  &&
  =
  \frac{1}{t_n} \sum_{t=0}^{t_n-1}
  (  F_n(\bw^{(t)}) -
  F_n(\bw^*))
  \\
  &&
  =
  \frac{1}{t_n} \sum_{t=0}^{t_n-1}
  (  F_{n,lin, \bw^{(t)}}(\bw^{(t)}) -
  F_{n,lin,\bw^{(t)}}(\bw^{*}))
  +
  \frac{1}{t_n} \sum_{t=0}^{t_n-1}
      ( F_{n,lin,\bw^{(t)}}(\bw^{*})-F_n(\bw^*))
      \\
      &&
      =: T_{1,n} + T_{2,n}.
    \end{eqnarray*}
    It is shown in the proof of Lemma 1 in Kohler (2024)
    that
    assumption (\ref{le1eq1}) implies
    \[
|    f_{\bw}(x)-f_{lin,\bw_0,\bw}(x)| \leq  \frac{1}{2} \cdot C_n \cdot \| \bw - \bw_0\|^2
\]
for all $x \in \{X_1, \dots, X_n\}$ and all $\bw_0, \bw \in A$.

Using (\ref{le1eq3})--(\ref{le1eq5}) we can conclude
\begin{eqnarray*}
  &&
|  F_n(\bw^*)
  -
    F_{n,lin,\bw^{(t)}}(\bw^*)
    |
    \\
    &&
    \leq
    \frac{1}{n} \sum_{i=1}^n |Y_i - f_{\bw^*}(X_i) + Y_i - f_{lin,\bw^{(t)},\bw^*}(X_i)|
    \cdot
    |f_{\bw^*}(X_i)  - f_{lin,\bw^{(t)},\bw^*}(X_i)|
    \\
    &&
    \leq
    \frac{1}{n} \sum_{i=1}^n (4 \cdot \beta_n +   
    \frac{1}{2} \cdot C_n \cdot \| \bw^* - \bw^{(t)}\|^2) \cdot
    \frac{1}{2} \cdot C_n \cdot \| \bw^* - \bw^{(t)}\|^2
    \\
    &&
    \leq
    \frac{1}{n} \sum_{i=1}^n (4 \cdot \beta_n + \frac{1}{2} \cdot  C_n \cdot
    4 \delta_n^2) \cdot \frac{1}{2} \cdot C_n  \| \bw^* - \bw^{(t)}\|^2
    \\
    &&
    \leq
    3 \cdot \beta_n \cdot C_n \cdot \| \bw^* - \bw^{(t)} \|^2.
  \end{eqnarray*}
This proves
\[
T_{2,n}
=
\frac{1}{t_n} \sum_{t=0}^{t_n-1}
      ( F_{n,lin,\bw^{(t)}}(\bw^{*})-F_n(\bw^*))
   \leq
    3 \cdot \beta_n \cdot C_n \cdot \frac{1}{t_n}  \sum_{t=0}^{t_n-1}
  \| \bw^* - \bw^{(t)} \|^2,
\]
hence it suffices to show
\begin{equation}
  \label{le1peq1}
T_{1,n}
\leq
\frac{\| \bw^* - \bw^{(0)} \|^2}{2 \cdot \lambda_n \cdot t_n}
+
\frac{1}{2} \cdot \lambda_n \cdot
\frac{1}{t_n} \sum_{t=0}^{t_n-1} \|\nabla_\bw F_n(\bw^{(t)}) \|^2
.
\end{equation}

    The convexity of $ F_{n,lin,\bw^{(t)}}( \bw )$ together with $\bw^* \in A$
    implies
    \begin{eqnarray*}
      &&
      F_{n,lin,\bw^{(t)}}(\bw^{(t)}) - F_{n,lin,\bw^{(t)}}(\bw^{*})
      \\
      &&
      \leq
\;      < \nabla_\bw       F_{n,lin,\bw^{(t)}}(\bw^{(t)}), \bw^{(t)}-\bw^* >
      \\
      &&
      =
  \;          < \nabla_\bw       F_{n}(\bw^{(t)}), \bw^{(t)}-\bw^* >
            \\
            &&
            =
            \frac{1}{2 \cdot \lambda_n}
            \cdot 2 \cdot < \lambda_n  \cdot \nabla_\bw    F_{n}(\bw^{(t)}), \bw^{(t)}-\bw^* >
            \\
            &&
            =
             \frac{1}{2 \cdot \lambda_n}
             \cdot
             \left(
\| \bw^{(t)}-\bw^*\|^2 - \|\bw^{(t)}-\bw^* - \lambda_n  \cdot \nabla_\bw    F_{n}(\bw^{(t)})\|^2 + \| \lambda_n  \cdot \nabla_\bw    F_{n}(\bw^{(t)})\|^2
\right)
\\
&&
     =
             \frac{1}{2 \cdot \lambda_n}
             \cdot
             \left(
\| \bw^{(t)}-\bw^*\|^2 - \|\bw^{(t)} - \lambda_n  \cdot \nabla_\bw    F_{n}(\bw^{(t)})- \bw^* \|^2 
\right)
+
\frac{1}{2} \cdot \lambda_n \cdot
\| \nabla_\bw    F_{n}(\bw^{(t)})\|^2
\\
&&
\leq
             \frac{1}{2 \cdot \lambda_n}
             \cdot
             \left(
             \| \bw^{(t)}-\bw^*\|^2 - \|Proj_A \left(
             \bw^{(t)} - \lambda_n  \cdot \nabla_\bw
             F_{n}(\bw^{(t)}) \right)- \bw^* \|^2 
\right)
\\
&&
\quad
+
\frac{1}{2} \cdot \lambda_n \cdot
\| \nabla_\bw    F_{n}(\bw^{(t)})\|^2
\\
&&
=      \frac{1}{2 \cdot \lambda_n}
             \cdot
             \left(
             \| \bw^{(t)}-\bw^*\|^2 -
             \| \bw^{(t+1)}-\bw^*\|^2 \right)
             +
\frac{1}{2} \cdot \lambda_n \cdot
\| \nabla_\bw    F_{n}(\bw^{(t)})\|^2.
      \end{eqnarray*}
    This implies
    \begin{eqnarray*}
      T_{1,n}
      & \leq &
      \frac{1}{t_n} \sum_{t=0}^{t_n-1}
      \left(
 \frac{1}{2 \cdot \lambda_n}
             \cdot
             \left(
             \| \bw^{(t)}-\bw^*\|^2 -
             \| \bw^{(t+1)}-\bw^*\|^2 \right)
             +
\frac{1}{2} \cdot \lambda_n \cdot
\| \nabla_\bw    F_{n}(\bw^{(t)})\|^2
\right)
\\
&\leq&
\frac{
\| \bw^{(0)}-\bw^*\|^2 
}{2 \cdot \lambda_n \cdot t_n}
    +
    \frac{1}{2} \cdot
    \frac{1}{t_n} \sum_{t=0}^{t_n-1}
    \lambda_n \cdot
\| \nabla_\bw    F_{n}(\bw^{(t)})\|^2,
      \end{eqnarray*}
which proves (\ref{le1peq1}). \hfill $\Box$%

\subsection{Proof of Theorem \ref{th1}}
We mimick the proof of Theorem 1 in Kohler (2024).

W.l.o.g. we assume
throughout the proof that $n$ is sufficiently large and that
$\|m\|_\infty \leq \beta_n$ holds.
Let $E_n$ be the event that
\[
\max_{i=1, \dots, n} |Y_i| \leq \sqrt{\beta_n}
\]
holds.

In the {\it first step of the proof} we show that on $E_n$
the conditions (\ref{se2eq6})--(\ref{se2eq8})
hold (provided we replace the constant $c_9$ in  (\ref{se2eq6})--(\ref{se2eq8})
by a larger constant, which we will denote again by $c_9$).

To show this it suffices to show that in case
\[
\hat{t}_n \geq n \cdot (\log n)^{c_{8}} \cdot K_n^3
\]
conditions  (\ref{se2eq6})--(\ref{se2eq8})
are satisfied. Observe that in this case we have
\[
n \cdot (\log n)^{c_{8}} \cdot K_n^3 \leq \hat{t}_n
\leq
2 \cdot n \cdot (\log n)^{c_{8}} \cdot K_n^3,
\]
which implies
\[
\lambda_n \cdot t_n
=
\frac{1}{\hat{t}_n} \cdot
\min \left\{
\hat{t}_n, \lceil (\log n)^{c_{8}} \cdot K_n^3 \rceil
\right\}
\geq
\frac{1}{ 2 \cdot n \cdot (\log n)^{c_{8}} \cdot K_n^3}
\cdot
\lceil (\log n)^{c_{8}} \cdot K_n^3 \rceil
\geq
\frac{1}{2 \cdot n}
\]
and
\[
\lambda_n \cdot t_n
=
\frac{1}{\hat{t}_n} \cdot
\min \left\{
\hat{t}_n, \lceil (\log n)^{c_{8}} \cdot K_n^3 \rceil
\right\}
\leq
\frac{
\lceil (\log n)^{c_{8}} \cdot K_n^3 \rceil
}{
n \cdot (\log n)^{c_{8}} \cdot K_n^3
}
\leq \frac{2}{n}.
\]

On $E_n$ we have
\[
F_n(\bw^{(0)})
=
\frac{1}{n} \sum_{i=1}^n |Y_i -0|^2 \leq \beta_n,
\]
hence
\[
\sqrt{
  8 \cdot \frac{t_n}{\hat{t}_n} \cdot
  \max \left\{
F_n(\bw^{(0)}),1
  \right\}
}
\leq 4 \cdot \frac{\sqrt{\beta_n}}{\sqrt{n}} \leq 1
\]
holds.

From this, $c_8 > 2L$ and the initial choice of $\bw^{(0)}$
we can conclude from Lemma 3 in Kohler (2024)
(which we apply with $\gamma_n^*=1$ and $B_n = c_{13} \cdot \log n +1$)
that
\[
\| \bw-\bw^{(0)}\|
\leq
\sqrt{
  2 \cdot \frac{t_n}{\hat{t}_n}
  \cdot \max \{
F_n(\bw^{(0)}),1
  \}
  }
\]
implies
\[
\| \nabla_\bw F_n(\bw) \| \leq
c_{18} \cdot (\log n)^L \cdot K_n^{3/2}
\leq
\sqrt{
  2 \cdot t_n \cdot \hat{t}_n \cdot
   \max \{
F_n(\bw^{(0)}),1
  \}
  },
\]
and by Lemma 5 in Kohler (2024) we see that
\[
\| \bw_1 - \bw^{(0)} \| \leq
\sqrt{
  8 \cdot \frac{t_n}{\hat{t}_n}
  \cdot \max \{
F_n(\bw^{(0)}),1
  \}
  }
\]
and
\[
\| \bw_2 - \bw^{(0)} \| \leq
\sqrt{
  8 \cdot \frac{t_n}{\hat{t}_n}
  \cdot \max \{
F_n(\bw^{(0)}),1
  \}
  }
\]
imply
\[
\| \nabla_\bw F_n(\bw_1) - \nabla_\bw F_n(\bw_2) \|
\leq
c_{19} \cdot K_n^{3/2} \cdot (\log n)^{2L} \cdot \| \bw_1 - \bw_2 \|
\leq
\hat{t}_n \cdot  \| \bw_1 - \bw_2 \|.
\]
Hence the assumptions of Lemma 4 in
Kohler (2024) are satisfied, and from this lemma we immediately
get
\[
\| \bw^{(t)} - \bw^{(0)} \|
\leq
\sqrt{
  2 \cdot \frac{t_n}{\hat{t}_n}
  \cdot \max \{
F_n(\bw^{(0)}),1
  \}
  }
\leq
\sqrt{
  \frac{4 \cdot \beta_n}{n}
  }
\quad (t=1, \dots, t_n)
\]
and
\[
F_n (\bw^{(t)}) \leq F_n (\bw^{(t-1)})
\quad (t=1, \dots, t_n),
\]
which implies
(\ref{se2eq7}) and (\ref{se2eq8}).
Furthermore, another application of Lemma 3 in Kohler (2024) yields
\[
\frac{1}{t_n}
\cdot \sum_{t=0}^{t_n -1}
\lambda_n \cdot \left\|
\nabla_\bw F_n (\bw^{(t)})
\right\|^2
\leq
\frac{1}{\hat{t}_n} \cdot c_{20} \cdot (\log n)^{2L}\cdot K_n^{3}
\leq
\frac{c_9}{n}
\]
(where we have used $c_8 > 2L$),
which completes the first step of the proof.

In the {\it second step of the proof} we decompose the $L_2$ error
in a sum of several terms. To do this we set
$m_{\beta_n}(x)= \EXP\{ T_{\beta_n} Y | X=x\}$
and observe
\begin{eqnarray*}
&&
\int | m_n(x)-m(x)|^2 \PROB_X (dx)
\\
&&
=
\left(
\EXP \left\{ |m_n(X)-Y|^2 | \D_n \right\}
-
\EXP \{ |m(X)-Y|^2\}
\right)
\cdot 1_{E_n}
+
\int | m_n(x)-m(x)|^2 \PROB_X (dx)
\cdot 1_{E_n^c}
\\
&&
=
\Big[
\EXP \left\{ |m_n(X)-Y|^2 | \D_n \right\}
-
\EXP \{ |m(X)-Y|^2\}
\\
&&
\hspace*{2cm}
- \left(
\EXP \left\{ |m_n(X)-T_{\beta_n} Y|^2 | \D_n \right\}
-
\EXP \{ |m_{\beta_n}(X)- T_{\beta_n} Y|^2\}
\right)
\Big] \cdot 1_{E_n}
\\
&&
\quad +
\Big[
\EXP \left\{ |m_n(X)-T_{\beta_n} Y|^2| \D_n \right\}
-
\EXP \{ |m_{\beta_n}(X)- T_{\beta_n} Y|^2\}
\\
&&
\hspace*{2cm}
-
2 \cdot \frac{1}{n} \sum_{i=1}^n
\left(
|m_n(X_i)-T_{\beta_n} Y_i|^2
-
|m_{\beta_n}(X_i)- T_{\beta_n} Y_i|^2
\right)
\Big] \cdot 1_{E_n}
\\
&&
\quad
+\Big[
2 \cdot \frac{1}{n} \sum_{i=1}^n
|m_n(X_i)-T_{\beta_n} Y_i|^2
-
2 \cdot \frac{1}{n} \sum_{i=1}^n
|m_{\beta_n}(X_i)- T_{\beta_n} Y_i|^2
\\
&&
\hspace*{2cm}
- \left(
2 \cdot \frac{1}{n} \sum_{i=1}^n
|m_n(X_i)-Y_i|^2
-
2 \cdot \frac{1}{n} \sum_{i=1}^n
|m(X_i)- Y_i|^2
\right)
\Big] \cdot 1_{E_n}
\\
&&
\quad
+
\Big[
2 \cdot \frac{1}{n} \sum_{i=1}^n
|m_n(X_i)-Y_i|^2
-
2 \cdot \frac{1}{n} \sum_{i=1}^n
|m(X_i)- Y_i|^2
\Big] \cdot 1_{E_n}
\\
&&
\quad
+
\int | m_n(x)-m(x)|^2 \PROB_X (dx)
\cdot 1_{E_n^c}
\\
&&
=: \sum_{j=1}^5 T_{j,n}.
\end{eqnarray*}
In the remainder of the proof we bound
\[
\EXP T_{j,n}
\]
for $j \in \{1, \dots, 5\}$.

In the {\it third step of the proof} we show
\[
\EXP T_{j,n} \leq c_{21} \cdot \frac{\log n}{n} \quad
\mbox{for } j \in \{1,3\}.
\]
This follows as in the proof of Lemma 1 in Bauer and Kohler (2019).

In the {\it fourth step of the proof} we show
\[
\EXP T_{5,n} \leq c_{22} \cdot \frac{(\log n)^2}{n}.
\]

The definition of $m_n$ implies $\int |m_n(x)-m(x)|^2 \PROB_X (dx) \leq
4 \cdot c_{12}^2 \cdot (\log
n)^2$, hence
\begin{eqnarray}
\PROB(E_n^c)
&\leq&
\PROB\{ \max_{i=1, \dots, n} |Y_i| > \sqrt{\beta_n}
\}
 \leq
n \cdot\PROB\{  |Y| > \sqrt{\beta_n}
\}
\nonumber
\\
& \leq &
n \cdot
\frac{\EXP\{ \exp(c_7 \cdot Y^2)}{ \exp( c_7 \cdot \beta_n)}
\leq
\frac{c_{23}}{n^2}
\label{pth1eq1}
\end{eqnarray}
where the last inequality holds because of (\ref{th1eq1})  and
$c_7 \cdot c_{12} \geq 3$, implies the assertion.

Let $\epsilon >0$ be arbitrary.
In the {\it fifth step of the proof} we show
\[
\EXP T_{2,n} \leq
c_{24} \cdot
\frac{ A_n^d \cdot B_n^{(L-1) \cdot d}}{n^{1-\epsilon}}
.
\]
Let $\W_n$ be the set of all weight vectors
$(w_{i,j,k}^{(l)})_{i,j,k,l}$ which satisfy
\[
| w_{k,1,1}^{(L)}| \leq c_{25} \quad (k=1, \dots, K_n),
\]
\[
|w_{k,i,j}^{(l)}| \leq c_{26} \cdot B_n \quad (l=1, \dots, L-1)
\]
and
\[
|w_{k,i,j}^{(0)}| \leq c_{27} \cdot A_n.
\]
By the first step of the proof we know that on $E_n$
condition
(\ref{se2eq8}) holds.
From this and the initial
choice of $\bw^{(0)}$ we can 
conclude that on $E_n$ we have
\[
\bw^{(t_n)} \in \W_n.
\]
Hence, for any $u>0$ we get
\begin{eqnarray*}
&&
\PROB \{ T_{2,n} > u \}
\\
&&
\leq
\PROB \Bigg\{
\exists f \in \F_n :
\EXP \left(
\left|
\frac{f(X)}{\beta_n} - \frac{T_{\beta_n}Y}{\beta_n}
\right|^2
\right)
-
\EXP \left(
\left|
\frac{m_{\beta_n}(X)}{\beta_n} - \frac{T_{\beta_n}Y}{\beta_n}
\right|^2
\right)
\\
&&\hspace*{3cm}-
\frac{1}{n} \sum_{i=1}^n
\left(
\left|
\frac{f(X_i)}{\beta_n} - \frac{T_{\beta_n}Y_i}{\beta_n}
\right|^2
-
\left|
\frac{m_{\beta_n}(X_i)}{\beta_n} - \frac{T_{\beta_n}Y_i}{\beta_n}
\right|^2
\right)
\Bigg\}
\\
&&\hspace*{2cm}
> \frac{1}{2} \cdot
\left(
\frac{u}{\beta_n^2}
+
\EXP \left(
\left|
\frac{f(X)}{\beta_n} - \frac{T_{\beta_n}Y}{\beta_n}
\right|^2
\right)
-
\EXP \left(
\left|
\frac{m_{\beta_n}(X)}{\beta_n} - \frac{T_{\beta_n}Y}{\beta_n}
\right|^2
\right)
\right) \Bigg\},
\end{eqnarray*}
where
\[
\F_n = \left\{ T_{\beta_n} f_\bw \quad : \quad \bw \in \W_n \right\}.
\]
By Lemma 12 in Kohler (2024) we get 
\begin{eqnarray*}
&&
\Nu_1 \left(
\delta , \left\{
\frac{1}{\beta_n} \cdot f : f \in \F_n
\right\}
, x_1^n
\right)
\leq
\Nu_1 \left(
\delta \cdot \beta_n , \F_n
, x_1^n
\right)
\\
&&
\leq
\left(
\frac{ c_{28}}{\delta}
\right)^{
c_{29} \cdot   A_n^{d} \cdot B_n^{(L-1) \cdot d}
\cdot
   \left(\frac{K_n \cdot c_{30}}{\beta_n \cdot \delta}\right)^{d/k} + c_{31}
  }.
\end{eqnarray*}
By choosing $k$ large enough we get for $\delta>1/n^2$
\[
\Nu_1 \left(
\delta , \left\{
\frac{1}{\beta_n} \cdot f : f \in \F_n
\right\}
, x_1^n
\right)
\leq
c_{32} \cdot n^{ c_{33} \cdot A_n^d \cdot  B_n^{(L-1) \cdot d} \cdot n^{\epsilon/2}}.
\]
This together with Theorem 11.4 in Gy\"orfi et al. (2002) leads for $u
\geq 1/n$ to
\[
\PROB\{T_{2,n}>u\}
\leq
14 \cdot
c_{32} \cdot n^{ c_{33} \cdot A_n^d \cdot  B_n^{(L-1) \cdot d} \cdot n^{\epsilon/2}}
\cdot
\exp \left(
- \frac{n}{5136 \cdot \beta_n^2} \cdot u
\right).
\]
For $\epsilon_n \geq 1/n$ we can conclude
\begin{eqnarray*}
\EXP \{ T_{2,n} \}
& \leq &
\epsilon_n + \int_{\epsilon_n}^\infty \PROB\{ T_{2,n}>u \} \, du
\\
& \leq &
\epsilon_n
+
14 \cdot
c_{32} \cdot n^{ c_{33} \cdot A_n^d \cdot  B_n^{(L-1) \cdot d} \cdot n^{\epsilon/2}}
\cdot
\exp \left(
- \frac{n}{5136 \cdot \beta_n^2} \cdot \epsilon_n
\right)
\cdot
\frac{5136 \cdot \beta_n^2}{n}.
\end{eqnarray*}
Setting
\[
\epsilon_n = \frac{5136 \cdot \beta_n^2}{n}
\cdot
c_{33} \cdot A_n^d \cdot  B_n^{(L-1) \cdot d} \cdot n^{\epsilon/2}
\cdot \log n
=
\frac{5136 \cdot \beta_n^2}{n}
\cdot
\log
\left(
n^{
c_{33} \cdot A_n^d \cdot  B_n^{(L-1) \cdot d} \cdot n^{\epsilon/2}
  }
\right)
\]
yields the assertion of the fourth step of the proof.

In the {\it sixth step of the proof} we show
\begin{eqnarray*}
&&
  \EXP \{ T_{4,n} \} 
  \leq
c_{34} \cdot \Bigg(
\EXP \left\{
\inf_{\bw: \|\bw-\bw^{(0)}\| \leq \frac{1}{n}}
\int | f_\bw (x)-m(x)|^2 \PROB_X (dx)
\right\}
+
\frac{(\log n)^{2L+2}}{n}
\Bigg)
.
  \end{eqnarray*}
Using
\[
|T_{\beta_n} z - y| \leq |z-y|
\quad \mbox{for } |y| \leq \beta_n
\]
we get
\begin{eqnarray*}
&&
 T_{4,n}/2
\\
&&
=
\Big[ \frac{1}{n} \sum_{i=1}^n
|m_n(X_i)-Y_i|^2
-
 \frac{1}{n} \sum_{i=1}^n
|m(X_i)- Y_i|^2
\Big] \cdot 1_{E_n}
\\
&&
\leq
\Big[
\frac{1}{n} \sum_{i=1}^n
|f_{\bw^{(t_n)}}(X_i)-Y_i|^2
-
 \frac{1}{n} \sum_{i=1}^n
|m(X_i)- Y_i|^2
\Big] \cdot 1_{E_n}
\\
&&
\leq
\big[ F_n(\bw^{(t_n)})
-
 \frac{1}{n} \sum_{i=1}^n
|m(X_i)- Y_i|^2
\Big] \cdot 1_{E_n}.
\end{eqnarray*}
By the first step of the proof we know
that on $E_n$
\[
\bw^{(t)} \in A= \left\{ \bw \in \R^{J_n} \, : \, \| \bw - \bw^{(0)}\| \leq
\frac{c_{35} \cdot
  \sqrt{\log n}}{\sqrt{n}} \right\}
\]
holds for $t=1, \dots, t_n$.
If $\bw_1$ and $\bw_2$ satisfy
\[
\| \bw_i - \bw^{(0)} \|
\leq
\frac{
c_{36} \cdot \sqrt{\log n}
}{
\sqrt{n}
}
\quad (i \in \{1,2\}),
\]
then the initialization of $\bw^{(0)}$ implies
\[
| (\bw_i)_{k,1,1}^{(L)}| \leq c_{37}
\quad \mbox{and} \quad
| (\bw_i)_{k,i,j}^{(l)}| \leq c_{38} \cdot \log n \quad (l=1, \dots, L-1)
\]
for $i \in \{1,2\}$, and by
Lemma 2 in Kohler (2024) we can conclude
\[
  \sum_{j=1}^{J_n}
  \left|
  \frac{\partial}{\partial w^{(j)}}
  f_{\bw_1}(x)
-
\frac{\partial}{\partial w^{(j)}}
f_{\bw_2}(x)
\right|^2
   \leq c_{39} \cdot (\log n)^{4L} \cdot
  \| \bw_1 - \bw_2 \|^2
  \]
  for $x \in supp(X)$.
Application of Lemma \ref{le1} with $A$ defined as above 
and $C_n = c_{40} \cdot (\log n)^{2L}$
yields because of $\lambda_n \cdot t_n \geq 1/2n$
(which follows from the first step of the proof)
\begin{eqnarray*}
&&
 T_{4,n}/2
\\
&&
\leq
\Big[
\frac{1}{n} \sum_{i=1}^n
|f_{\bw^*}(X_i)-Y_i|^2
+
c_{41} \cdot \frac{(\log n)^{2L+2}}{n}
-
 \frac{1}{n} \sum_{i=1}^n
|m(X_i)- Y_i|^2
\Big] \cdot 1_{E_n}
\\
&&
\leq
\frac{1}{n} \sum_{i=1}^n
|f_{\bw^*}(X_i)-Y_i|^2
-
 \frac{1}{n} \sum_{i=1}^n
 |m(X_i)- Y_i|^2
 +
 c_{41} \cdot \frac{(\log n)^{2L+2}}{n}
 \\
 &&
 \quad
 +
  \frac{1}{n} \sum_{i=1}^n
|m(X_i)- Y_i|^2
 \cdot 1_{E_n^c}
\end{eqnarray*}
for any $\bw^*$ with $\|\bw^*-\bw^{(0)}\| \leq 1/n$.
Hence using (\ref{pth1eq1}) we can conclude
\[
\EXP\{ T_{4,n}/2 | \bw^{(0)} \}
\leq
\int | f_{\bw^*}(x)-m(x)|^2 \PROB_X (dx) +
 c_{42} \cdot \frac{(\log n)^{2L+2}}{n}
\]
for any $\bw^*$ with $\|\bw^*-\bw^{(0)}\| \leq 1/n$, which implies
\[
\EXP\{ T_{4,n}/2 | \bw^{(0)} \}
\leq
\inf_{\bw : \| \bw-\bw^{(0)}\| \leq 1/n}
\int | f_{\bw}(x)-m(x)|^2 \PROB_X (dx) +
 c_{42} \cdot \frac{(\log n)^{2L+2}}{n}
\]
and
\[
\EXP\{ T_{4,n}/2 \}
\leq
\EXP \left\{
\inf_{\bw : \| \bw-\bw^{(0)}\| \leq 1/n}
\int | f_{\bw}(x)-m(x)|^2 \PROB_X (dx)
\right\}
+
 c_{42} \cdot \frac{(\log n)^{2L+2}}{n}.
\]

\hfill $\Box$

\subsection{An auxiliary result for the proof of Corollary \ref{co1}}

In the proof of Corollary \ref{co1} we will need the following
result concerning the approximation of $(p,C)$--smooth functions
by neural networks with bounded weights. 

\begin{lemma}
  \label{le2}
Let $d \in \N$,  $p=q+\beta$ where
$\beta \in (0,1]$ and $q \in \N_0$, $C>0$,
$A \geq 1$
and
$A_n, B_n, \gamma_n^* \geq 1$.
For $L,r,K \in \N$
let $\F$ be the set of all networks $f_{\bw}$ defined by
(\ref{se2eq1})--(\ref{se2eq3}) with $K_n$ replaced by $r$, where
the weight vector satisfies
\[
|w_{i,j}^{(0)}| \leq A_n, \quad
|w_{i,j}^{(l)}| \leq B_n \quad \mbox{and} \quad
|w_{i,j}^{(L)}| \leq \gamma_n^*
\]
for all $l \in \{1, \dots, L-1\}$ and all $i,j$, and set
\[
\HH = \left\{ \sum_{k=1}^{K^d} f_k \quad : \quad f_k \in \F \quad (k=1, \dots, K)
\right\}.
\]
Let $L,r \in \N$ with
\[
L \geq \lceil \log_2(q+d) \rceil
\quad
\mbox{and}
\quad
r \geq 2 \cdot (2p+d) \cdot (q+d),
\]
and set
\[
 A_n =A \cdot K \cdot \log K, \quad B_n=c_{43}
\quad \mbox{and} \quad 
\gamma_n^*=c_{44} \cdot K^{q+d}.
\]
Assume $K \geq c_{45}$ for $c_{45}>0$ sufficiently large.
Then there exists for any $(p,C)$--smooth $f:\R^d \rightarrow \R$
a neural network $h \in \HH$ such that
\[
\sup_{x \in [-A,A)^d} |f(x)-h(x)|
\leq
\frac{c_{46}}{K^{p}}.
\]
\end{lemma}

\noindent
{\bf Proof.} See Theorem 3 in Kohler (2024). \hfill $\Box$

\subsection{Proof of Corollary \ref{co1}}
In the proof we will use arguments from the proof of Theorem 1
in Kohler (2024).

W.l.o.g. we assume
throughout the proof that $n$ is sufficiently large and that
$\|m\|_\infty \leq \beta_n$ holds.
Let $A>0$ with $supp(X) \subseteq [-A,A]^d$.
Set
\[
\tilde{K}_n= \left\lceil c_{47} \cdot n^{\frac{d}{2p+d}} \right\rceil 
\]
and
\[
N_n= \left\lceil c_{48} \cdot n^{4+\frac{d}{2p+d}}  \right\rceil
\]
and let $\bw^*$ be a weight vector of a neural networks
where the results of $N_n \cdot \tilde{K}_n \cdot r$ in parallel computed neural
networks with $L$ hidden layers and $r$ neurons per layer are
computed such that the corresponding network
\[
f_{\bw^*}(x)= \sum_{k=1}^{N_n \cdot \tilde{K}_n \cdot r} (\bw^*)_{k,1,1}^{(L)} \cdot f_{\bw^*,k,1}^{(L)}(x)
\]
satisfies
\begin{equation}
\label{pco1eq1}
\sup_{x \in [-A,A]^d} |f_{\bw^*}(x)-m(x)| \leq \frac{c_{49}}{\tilde{K}_n^{p/d} }
\end{equation}
and
\[
|(\bw^*)_{k,1,1}^{(L)}| \leq \frac{c_{50} \cdot  \tilde{K}_n^{(q+d)/d}}{N_n}
\quad (k=1, \dots, N_n \cdot \tilde{K}_n \cdot r).
\]
Note that such a network exists according to Lemma \ref{le2}
if we repeat in the outer sum of the function space $\HH$
each of the $f_k$'s in Lemma \ref{le2} $N_n$--times
with outer weights divided by $N_n$.
Set
\[
\epsilon_n = \frac{c_{51}}{ n \cdot \sqrt{N_n \cdot \tilde{K}_n}}
\geq
\frac{c_{52}}{n^{4}}.
\]

  Let $E_n$ be the event that the weight vector $\bw^{(0)}$
            satisfies
            \[
            | (\bw^{(0)})_{j_s,k,i}^{(l)}-(\bw^*)_{s,k,i}^{(l)}| \leq \epsilon_n
            \quad \mbox{for all } l \in \{0, \dots, L-1\},
            s \in \{1, \dots, N_n \cdot \tilde{K}_n \cdot r\}
            \]
            for some pairwise distinct $j_1, \dots, j_{N_n \cdot \tilde{K}_n \cdot r}
            \in \{1, \dots, K_n\}$.

            In the {\it first step of the proof} we show
            \begin{equation}
              \label{pco1eq2}
              \PROB(E_n^c) \leq c_{53} \cdot n^{6} \cdot \exp( - n^{0.5}).
              \end{equation}
To do this, we consider a sequential choice of the weights of 
$K_n$ fully connected neural networks. The probability that the weights in
the first of these networks differ in all components at most by $\epsilon_n$
from $((\bw^*)_{1,i,j}^{(l)})_{i,j,l: l<L}$  is
for large $n$ bounded from below by
\begin{eqnarray*}
  &&
  \left( \frac{c_{52}}{2 \cdot c_{54} \cdot (\log n)  \cdot n^{4}}
\right)^{r \cdot (r+1) \cdot (L-1)}
\cdot
\left(
\frac{ c_{52} }{2 \cdot c_{55} \cdot (\log n) \cdot n^{1/(2p+d)} \cdot n^{4}}
\right)^{r \cdot (d+1)}
\\
&&
\geq
n^{-r \cdot (r+1) \cdot (L-1) \cdot 4- r \cdot 4 \cdot (d+1) - r \cdot (d+1)  /(2p+d) - 0.5}.
\end{eqnarray*}
Hence probability that none of the first $n^{
r \cdot (r+1) \cdot (L-1) \cdot 4+ r \cdot 4 \cdot (d+1)+ r \cdot (d+1) /(2p+d) +1}$ neural networks satisfies this condition is for large $n$
 bounded above by
 \begin{eqnarray*}
   &&
(1 -  n^{-r \cdot (r+1) \cdot (L-1) \cdot 4- r \cdot 4 \cdot (d+1)- r \cdot (d+1) /(2p+d) -0.5}) ^{n^{r \cdot (r+1) \cdot (L-1) \cdot 4 + r \cdot 4 \cdot (d+1) + r \cdot (d+1) / (2p+d) +1}}\\
&&\leq
\left(\exp \left(
-  n^{-r \cdot (r+1) \cdot (L-1) \cdot 4- r \cdot 4 \cdot (d+1) - r \cdot (d+1) / (2p+d) -0.5}
\right)
\right) ^{n^{r \cdot (r+1) \cdot (L-1) \cdot 4+ r \cdot 4 \cdot (d+1) + r \cdot (d+1) / (2p+d) +1}}
\\
&&=
\exp( -  n^{0.5}).
\end{eqnarray*}
 Since we have $K_n \geq n^{r \cdot (r+1) \cdot (L-1) \cdot 4+ r \cdot 4 \cdot (d+1) +
   r \cdot (d+1) / (2p+d) +1} \cdot N_n \cdot \tilde{K}_n \cdot r$
 for $n$ large we can successively
use the same construction for all of $N_n \cdot \tilde{K}_n \cdot r$ weights and we can conclude:
The probability that there exists $k \in \{1, \dots, N_n \cdot \tilde{K}_n \cdot r\}$
such that
none of the $K_n$ weight vectors of the fully
connected neural network differs by at most $\epsilon_n$ from
$((\bw^*)_{k,i,j}^{(l)})_{i,j,l:l<L}$ is for large $n$ bounded from above by
\begin{eqnarray*}
&&
N_n \cdot \tilde{K}_n \cdot r \cdot \exp( -  n^{0.5})
\leq  c_{56} \cdot n^{6} \cdot \exp( - n^{0.5}),
\end{eqnarray*}
which implies the assertion of the first step of the proof.

In the {\it second step of the proof} we show
\begin{equation}
  \label{pco1eq3}
  \EXP \left\{
\inf_{\bw: \|\bw-\bw^{(0)}\| \leq \frac{1}{n}}
\int | f_\bw (x)-m(x)|^2 \PROB_X (dx)
\right\}
\leq
c_{57} \cdot n^{- \frac{2p}{2p+d}}.
\end{equation}

On $E_n$ we have
\begin{eqnarray*}
  \| \bw^* - \bw^{(0)} \|^2
  & \leq &
  \sum_{k=1}^{N_n \cdot \tilde{K}_n \cdot r}
  | (\bw^*)_{k,1,1}^{(L)}|^2
  +
  N_n \cdot \tilde{K}_n \cdot r \cdot
  L \cdot (r+d)^2 \cdot \epsilon_n^2
  \\
  & \leq &
  \frac{c_{50}^2 \cdot \tilde{K}_n^{1 + 2 \cdot \frac{p+d}{d}}}{N_n}
  +
  \frac{c_{51}^2 \cdot r \cdot L \cdot (r+d)^2 }{n^2}
  \\
  & \leq & \frac{1}{n^2},
\end{eqnarray*}
provided $n$ is sufficiently large.
This implies
\begin{eqnarray*}
  &&
\EXP \left\{
\inf_{\bw: \|\bw-\bw^{(0)}\| \leq \frac{1}{n}}
\int | f_\bw (x)-m(x)|^2 \PROB_X (dx)
\right\}
\\
&&
\leq
\EXP \left\{
\int | f_{\bw^*} (x)-m(x)|^2 \PROB_X (dx) \cdot 1_{E_n}
\right\}
\\
&&
\quad
+
\EXP \left\{
\inf_{\bw: \|\bw-\bw^{(0)}\| \leq \frac{1}{n}}
\int | f_\bw (x)-m(x)|^2 \PROB_X (dx) \cdot 1_{E_n^c}
\right\}
\\
&&
\leq
\int | f_{\bw^*} (x)-m(x)|^2 \PROB_X (dx)
+
\int | 0-m(x)|^2 \PROB_X (dx) \cdot \PROB\{E_n^c\}
\\
&&
\leq
\int | f_{\bw^*} (x)-m(x)|^2 \PROB_X (dx)
+ \frac{c_{58}}{n},
\end{eqnarray*}
where the second last inequality followed from $f_{\bw^{(0)}}(x)=0$ for all
$x \in \R^d$.
Application of (\ref{pco1eq1}) yields the assertion.

In the {\it third step of the proof} we show the assertion.

Application of Theorem \ref{th1} with $\epsilon$ replaced by $\epsilon/2$
together with the result of the
second step of the proof yields
\begin{eqnarray*}
  &&
  \EXP \int | m_n(x)-m(x)|^2 \PROB_X (dx)
  \\
  &&
\leq
c_{59} \cdot \Bigg(
\EXP \left\{
\inf_{\bw: \|\bw-\bw^{(0)}\| \leq \frac{1}{n}}
\int | f_\bw (x)-m(x)|^2 \PROB_X (dx)
\right\}
+
\frac{A_n^d \cdot B_n^{(L-1) \cdot d}}{n^{1-\epsilon/2}}
\Bigg)
\\
&&
\leq
c_{60} \cdot \left(
n^{- \frac{2p}{2p+d}} +
\frac{A_n^d \cdot B_n^{(L-1) \cdot d}}{n^{1-\epsilon/2}}
\right)
\\
&&
\leq
c_{61} \cdot \left(
n^{- \frac{2p}{2p+d}} +
\frac{ n^{\frac{d}{2p+d}} \cdot (\log n)^{L \cdot d}}{n^{1-\epsilon/2}}
\right)
\\
&&
\leq
c_{62} \cdot 
n^{- \frac{2p}{2p+d} + \epsilon}.
\end{eqnarray*}

\hfill $\Box$

\end{document}